\documentclass[a4paper,11pt]{amsart}
\usepackage{amsmath,amsthm,amssymb,mathtools}
\usepackage[mathscr]{eucal}
 \usepackage{cite}
\usepackage{upgreek}
\usepackage[bookmarks,bookmarksnumbered]{hyperref}
\usepackage{enumerate}
\usepackage{amsmath}
\usepackage{mathrsfs}
\usepackage{blindtext}
\usepackage{scrextend}
\usepackage{enumitem}
\addtokomafont{labelinglabel}{\sffamily}
\usepackage{color}
\usepackage[at]{easylist}

\numberwithin{equation}{section}

\theoremstyle{plain}
\newtheorem{main}{Theorem}
\newtheorem{mcor}[main]{Corollary}
\newtheorem{mprop}[main]{Proposition}

\newtheorem{theorem}{Theorem}[section]
\newtheorem{claim}[theorem]{Claim}
\newtheorem{lemma}[theorem]{Lemma}
\newtheorem{proposition}[theorem]{Proposition}
\newtheorem{corollary}[theorem]{Corollary}
\theoremstyle{definition}
\newtheorem{definition}[theorem]{Definition}
\newtheorem*{definition*}{Definition}
\newtheorem{example}[theorem]{Example}

\newtheorem{remark}[theorem]{Remark}

\newcommand{\N}{\mathbb{N}}
\newcommand{\R}{\mathbb{R}}
\newcommand{\C}{\mathbb{C}}
\newcommand{\Z}{\mathbb{Z}}

\newcommand{\bB}{\mathbb{B}}

\newcommand{\cH}{\mathcal{H}}

\newcommand{\cK}{\mathcal{K}}

\newcommand{\cU}{\mathcal{U}}

\newcommand{\bcH}{\overline{\mathcal{H}}}
\newcommand{\ot}{\otimes}
\newcommand{\otb}{\bar{\otimes}}

\newcommand{\SO}{\operatorname{SO}}
\newcommand{\SU}{\operatorname{SU}}

\newcommand{\Ker}{\operatorname{Ker}}

\newcommand{\eps}{\varepsilon}

\newcommand{\norm}[1]{\left\|#1\right\|}

\begin{document}

\title[Cohomological obstructions to lifting properties for full group C$^*$-algebras]
{Cohomological obstructions to lifting properties \\ for  full  C$^*$-algebras of property (T) groups}

\author[A. Ioana]{Adrian Ioana}
\address{Department of Mathematics, University of California San Diego, 9500 Gilman Drive, La Jolla, CA 92093, USA}
\email{aioana@ucsd.edu}

\author[P. Spaas]{Pieter Spaas}
\address{Department of Mathematics, University of California Los Angeles, 520 Portola Plaza, Los Angeles, CA 90095, USA}
\email{pspaas@math.ucla.edu}

\author[M. Wiersma]{Matthew Wiersma}
\address{Department of Mathematics, University of California San Diego, 9500 Gilman Drive, La Jolla, CA 92093, USA}
\email{mtwiersma@ucsd.edu}

\thanks{A.I. was supported in part by NSF Career Grant DMS \#1253402 and NSF FRG Grant \#1854074.}
\begin{abstract} 
We develop a new method, based on non-vanishing of second cohomology groups, for proving the failure of lifting properties for full C$^*$-algebras  of  countable groups with (relative) property (T). We derive that the full   C$^*$-algebras of the groups $\mathbb Z^2\times\text{SL}_2(\mathbb Z)$ and $\text{SL}_n(\mathbb Z)$, for $n\geq 3$, do not have the local lifting property (\text{LLP}). We also prove that the full C$^*$-algebras of a large class of groups $\Gamma$ with property (T), including those such that $\text{H}^2(\Gamma,\mathbb R)\not=0$ or $\text{H}^2(\Gamma,\mathbb Z\Gamma)\not=0$, do not have the lifting property (\text{LP}). More generally, we show that the same holds if $\Gamma$ admits a probability measure preserving action with non-vanishing second $\mathbb R$-valued cohomology. Finally, we prove that the full C$^*$-algebra of any non-finitely presented property (T) group fails the LP.
\end{abstract}

\maketitle

\section{Introduction and statement of main results} 
The {\it local lifting property} (LLP) was introduced by Kirchberg in his landmark paper \cite{Ki93},
as a weaker, local version of the (global) 
{\it lifting property} (LP). For a unital C$^*$-algebra $A$, the LP requires that any unital completely positive (u.c.p.) map from $A$ into a quotient C$^*$-algebra admits a u.c.p. lift. We refer the reader to Section \ref{ssec:consequence} for the precise definitions. Over the years, both the LP and the LLP have proven to be very useful properties, for instance in connecting the Connes-Kirchberg problem with Tsirelson's problem (see, e.g., \cite[Chapter 16]{Pi20}).

The Choi-Effros lifting theorem \cite{CE76} shows that nuclear C$^*$-algebras, and hence C$^*$-algebras of amenable groups, have the LP. In \cite{Ki94b} Kirchberg showed that full C$^*$-algebras of countable free groups have the LP. 
Furthermore, in \cite{Ki93}, Kirchberg proved tensorial characterizations of the LLP and of Lance's weak expectation property (WEP), establishing a tensorial duality between 
 the LLP and the WEP. 
 This allowed him to prove equivalence of several C$^*$-algebraic statements,  including whether LLP $\Rightarrow$ WEP, and Connes' embedding problem. 
In the same paper, it is also proven that the LLP is preserved under various operations, such as extensions and taking tensor products with nuclear C$^*$-algebras. 
The LLP was then shown to be preserved under full free products in \cite{Pi96}.



Besides C$^*$-algebras obtained through the aforementioned constructions, only a few other examples of C$^*$-algebras with the LLP are known.
Only recently,  answering a question dating back to \cite{Ki93}, Pisier constructed the first example of a non-nuclear C$^*$-algebra with both LLP and WEP in \cite{Pi19}. Most recently, two more  concrete families of C$^*$-algebras with the LLP were constructed in \cite{Co20}.


On the other hand, a fundamental result of Junge and Pisier in \cite{JP95} shows that $\bB(\ell^2)$ does not have the LLP.  
Aside from this, not many examples of 
C$^*$-algebras that fail the LLP or the LP are known.
(Note that the difference between the two properties is in fact quite subtle: it is open whether there exist separable C$^*$-algebras with the LLP but not the LP, see \cite{Oz04b}.)
 In particular, little is known about which full group C$^*$-algebras do not have these properties. 
In \cite{Oz04a}, Ozawa showed the existence of full group C$^*$-algebras without the LP. A few years later, Thom  \cite{Th10} constructed two Connes-embeddable non-residually finite  property (T) groups whose full C$^*$-algebras fail the LLP.  The latter result gives the first, and so far only, examples of countable groups whose full C$^*$-algebras do not have the LLP.
Additionally, it provides the only concrete examples of full group C$^*$-algebras without the LP.  
In fact, the combination of the Connes embedding property, property (T) and non-residual finiteness was up to date the only known obstruction to the (L)LP for full group C$^*$-algebras.
This is a surprising situation because most full group C$^*$-algebras are expected to not have the LLP, as already mentioned by Ozawa in \cite{Oz04b}. More recently, during the workshop ``Amenability, coarse embeddability and fixed point properties" at MSRI in December 2016, Pisier explicitly
stated the problem of finding more full group C$^*$-algebras either with or without the LLP.

We make progress on this problem here by developing a new method to refute the (L)LP for full group C$^*$-algebras. This new approach is based on the establishment of certain cohomological obstructions to the LP and the LLP for groups with (relative) property (T). In particular, it allows us to provide many new, natural examples of groups whose full C$^*$-algebras fail the LLP and large additional classes of groups whose full C$^*$-algebras fail the LP. 
The following is our first main result.

\begin{main}\label{A} 
Let $\Gamma$ be a countable group and $\Lambda$ be a subgroup such that the pair $(\Gamma,\Lambda)$ has the relative property (T). 
Assume that there is a sequence of 2-cocycles $c_n\in\emph{Z}^2(\Gamma,\mathbb T)$, $n\in\mathbb N$, such that
\begin{enumerate}
\item the restriction of $c_n$ to $\Lambda$ is not a 2-coboundary, for every $n\in\mathbb N$, 
\item $\lim\limits_{n\rightarrow\infty}c_n(g,h)=1$, for every $g,h\in\Gamma$, and
\item for every $n\in\mathbb N$, there is a projective representation $\pi_n:\Gamma\rightarrow\mathcal U(\mathcal H_n)$ on a finite dimensional Hilbert space $\mathcal H_n$ such that $\pi_n(g)\pi_n(h)=c_n(g,h)\pi_n(gh)$, for every $g, h\in\Gamma$.

\end{enumerate}

Then C$^*(\Gamma)$ does not have the \emph{LLP}.
\end{main}

Given countable groups $\Lambda<\Gamma$,  the pair $(\Gamma,\Lambda)$ has the {\it relative property (T)} of Kazhdan-Margulis if any unitary representation of $\Gamma$ with almost invariant vectors has a non-zero $\Lambda$-invariant vector. If the pair $(\Gamma,\Gamma)$ has the relative property (T), then $\Gamma$ is said to have Kazhdan's {\it property (T)}.

The proof of Theorem \ref{A} combines a consequence of the \text{LLP} concerning ``almost homomorphisms" and a beautiful characterization of relative property (T) in terms of projective representations due to Nicoara, Popa and Sasyk \cite{NPS07}, see the comments at the end of the introduction.


The main example of a pair of groups to which Theorem \ref{A} applies is  $\Lambda=\mathbb Z^2<\Gamma=\mathbb Z^2\rtimes\text{SL}_2(\mathbb Z)$. The pair $(\Gamma,\Lambda)$ has the relative property (T) by \cite{Ka67,Ma82}, while the existence of a sequence of 2-cocycles satisfying the hypothesis of Theorem \ref{A} was pointed out in \cite{NPS07}. 
 Theorem \ref{A} therefore implies that $C^*(\Gamma)$ does not  have the LLP. More generally, using that the pair $(\mathbb Z^2\rtimes\Sigma,\mathbb Z^2)$ has the relative property (T) for any non-amenable subgroup $\Sigma<\text{SL}_2(\mathbb Z)$ \cite{Bu91}, 
that the pair $(R^2\rtimes\text{SL}_2(R), R^2)$ has the relative property (T) for any finitely generated commutative ring with unit $R$ \cite{Sh99}, 
 and that the LLP passes to subgroups (see Remark \ref{hereditary}), we derive the following:

\begin{mcor}\label{B} Let $\Gamma=\mathbb Z^2\rtimes\Sigma$, where $\Sigma<\emph{SL}_2(\mathbb Z)$ is a non-amenable subgroup. Then $C^*(\Gamma)$ does not have the \emph{LLP}.
Thus, $C^*(\mathbb Z^2\rtimes\emph{SL}_2(\mathbb Z))$ and $C^*(\emph{SL}_n(\mathbb Z))$ do not have the \emph{LLP}, for any $n\geq 3$. 

More generally, if $R$ is any finitely generated commutative ring with unit  such that $\{2x\mid x\in R\}$ is infinite, then $C^*(R^2\rtimes\emph{SL}_2(R))$ and $C^*(\emph{SL}_n(R))$ do not have the \emph{LLP}, for any $n\geq 3$.

Moreover, every countable group $\Delta$ that has a non-abelian free subgroup admits an action by automorphisms on a countable abelian group $A$ such that $C^*(A\rtimes\Delta)$ does not have the \emph{LLP}.
\end{mcor}

Corollary \ref{B} provides the first examples of residually finite groups (e.g., $\mathbb Z^2\rtimes\text{SL}_2(\mathbb Z)$ and $\text{SL}_n(\mathbb Z)$, $n\geq 3$) 
whose full C$^*$-algebras fail the LLP.  As observed in \cite{Oz04b}, for Connes-embeddable groups, the LLP implies Kirchberg's factorization property \cite{Ki94a}.
Since residually finite groups are Connes-embeddable and have the factorization property, Corollary \ref{B} implies that the LLP is strictly stronger than the factorization property for Connes-embeddable groups. 
Corollary \ref{B} also shows that  $C^*(\text{SL}_3(\mathbb Z))$ does not have the LLP and thus $\text{SL}_3(\mathbb Z)$ does not characterize the WEP (in the sense of \cite[Definition 3.4]{FKPT18}). This settles in the negative a question raised in \cite[page 114]{FKPT18}. Finally, Corollary \ref{B} also implies that the full C$^*$-algebras of the groups $\mathbb F_p[X]^2\rtimes\text{SL}_2(\mathbb F_p[X])$ and $\text{SL}_n(\mathbb F_p[X])$, for $n\geq 3$, do not have the LLP, for any prime $p\geq 3$. This allows to recover the examples of groups without the LLP exhibited in \cite[Section 2]{Th10}, except for the case $p=2$, and in \cite[Section 3]{Th10}, as these groups contain $\text{SL}_3(\mathbb  F_p[X])$ and $\text{SL}_3(\mathbb Z)$, respectively.

We now turn to results providing classes of full group C$^*$-algebras without the LP. 
If the hypothesis of Theorem \ref{A} is relaxed by removing assumption {\it (3)}, then by adapting the proof of Theorem \ref{A} we can show that $C^*(\Gamma)$ 
 fails the LP. Our next main result considerably generalizes this fact by allowing cocycles that arise from measure preserving actions. 
 
 Before stating this result in detail, we review some terminology.
Let $\Gamma$ be a countable group and $\mathcal A$ be an abelian group endowed with an action $\sigma:\Gamma\rightarrow\text{Aut}(\mathcal A)$. We denote
by $\text{Z}^2(\Gamma,\mathcal A)$ the group of $2$-cocycles, i.e., maps $c:\Gamma\times\Gamma\rightarrow\mathcal A$ satisfying $\sigma_g(c(h,k))c(g,hk)=c(g,h)c(gh,k)$, for all $g,h,k\in\Gamma$.  A 2-cocycle $c$ is a $2$-coboundary if there is a map $b:\Gamma\rightarrow\mathcal A$ such that $c(g,h)=b(g)\sigma_g(b(h))b(gh)^{-1}$, for all $g,h\in\Gamma$. We denote by $\text{B}^2(\Gamma,\mathcal A)$ the group of $2$-coboundaries and by $\text{H}^2(\Gamma,\mathcal A)=\text{Z}^2(\Gamma,\mathcal A)/\text{B}^2(\Gamma,\mathcal A)$ the second cohomology group of $\Gamma$ with coefficients in $\mathcal A$.

Let $(X,\mu)$ be a probability space, which will always be assumed standard, and $A$ be a Polish abelian group. Let $\text{L}^0(X,A)$ be the Polish group, with respect to the topology of convergence in measure, of equivalence classes of measurable functions $f:X\rightarrow A$, where two functions are equivalent if they coincide $\mu$-almost everywhere. Let $\Gamma\curvearrowright^{\sigma} (X,\mu)$ be a probability measure preserving (p.m.p.) action.  Then $\Gamma$ has a natural action on $\text{L}^0(X,A)$, which we still denote by $\sigma$ and is given by $\sigma_g(f)(x)=f(g^{-1}x)$.

\begin{main}\label{C}
Let $\Gamma$ be a countable group and $\Lambda$ be a subgroup such that the pair $(\Gamma,\Lambda)$ has the relative property (T).  Assume that there are a p.m.p. action $\Gamma\curvearrowright^{\sigma} (X,\mu)$ such that $\sigma_{|\Lambda}$ is ergodic, and  2-cocycles $c_n\in\emph{Z}^2(\Gamma,\emph{L}^0(X,\mathbb T))$ such that the restriction of $c_n$ to $\Lambda$ is not a 2-coboundary, for every $n\in\mathbb N$, and $\lim\limits_{n\rightarrow\infty}\|c_n(g,h)-1\|_{2}=0$, for every $g,h\in\Gamma$.

Then C$^*(\Gamma)$ does not have the \emph{LP}. 
Moreover, if  the twisted crossed product von Neumann algebra $\emph{L}^{\infty}(X)\rtimes_{\sigma,c_n}\Gamma$ embeds into $R^{\omega}$, for every $n\in\mathbb N$,
 then C$^*(\Gamma)$ does not have the \emph{LLP}.

\end{main}

Before presenting two applications of Theorem \ref{C}, let us make a remark on its hypothesis.

\begin{remark}\label{H^2(T)}
Theorem \ref{C} implies that if $\Gamma$ is a property (T) group such that $C^*(\Gamma)$ has the LP, then $\text{B}^2(\Gamma,\text{L}^0(X,\mathbb T))$ is an open subgroup of $\text{Z}^2(\Gamma,\text{L}^0(X,\mathbb T))$ and thus $\text{H}^2(\Gamma,\text{L}^0(X,\mathbb T))$ is countable, for every ergodic p.m.p. action $\Gamma\curvearrowright (X,\mu)$. 

It is a longstanding open problem, going back to Feldman and Moore's work \cite{FM77}, to calculate $\text{H}^2(\Gamma,\text{L}^0(X,\mathbb T))$ and, more generally, the higher cohomology groups $\text{H}^n(\Gamma,\text{L}^0(X,\mathbb T))$, for $n\geq 2$. It is known that $\text{H}^2(\Gamma,\text{L}^0(X,\mathbb T))=0$, for any free ergodic p.m.p. action $\Gamma\curvearrowright (X,\mu)$, if $\Gamma$ is an amenable group \cite{CFW81}, a free group or more generally a treeable group \cite{Ki17}. On the other hand, for free ergodic actions of property (T) groups, not a single calculation of $\text{H}^2$ is available (see \cite[6.6]{Po07a} and \cite[5.7.4]{AP18}).  In fact, no example is known of a free ergodic p.m.p. action $\Gamma\curvearrowright (X,\mu)$ of a property (T) group $\Gamma$ such that $\text{H}^2(\Gamma,\text{L}^0(X,\mathbb T))$ is trivial or even countable.

On the other hand, it is also a difficult task to produce examples of p.m.p. actions $\Gamma\curvearrowright (X,\mu)$ such that $\text{H}^2(\Gamma,\text{L}^0(X,\mathbb T))$ is uncountable or even non-trivial. Only a few years ago, Jiang \cite{Ji16} proved, by using Popa's cocycle superrigidity theorem \cite{Po07a},  that if $\Gamma$ is a property (T) group, $\lambda$ is the Haar measure of $\mathbb T$ and $\Gamma\curvearrowright (\mathbb T^\Gamma,\lambda^\Gamma)$ is the Bernoulli action, then $\text{H}^2(\Gamma,\mathbb T)\oplus\text{H}^2(\Gamma,\mathbb Z\Gamma)$ embeds into $\text{H}^2(\Gamma,\text{L}^0(\mathbb T^\Gamma,\mathbb T))$.
 In particular, if $\text{H}^2(\Gamma,\mathbb Z\Gamma)\not=0$, then $\text{H}^2(\Gamma,\text{L}^0(\mathbb T^\Gamma,\mathbb T))\not=0$.
\end{remark}

 By combining this result with Theorem \ref{C} and Popa's malleability property for Bernoulli actions, we derive the following:

\begin{mcor}\label{D}
Let $\Gamma$ be a countable group with property (T) such that $\emph{H}^2(\Gamma,\mathbb Z\Gamma)\not=0$.

Then $C^*(\Gamma)$ does not have the \emph{LP}.

\end{mcor}

Corollary \ref{D} implies that if $d\in(\frac{1}{3},\frac{1}{2})$, then for a random group $\Gamma$ in Gromov's density model at density $d$ (see \cite[Definition 7]{Ol05}), we have that
$C^*(\Gamma)$ does not have the LP, with overwhelming probability.
Indeed, as explained in the proof of \cite[Corollary 4.4]{Ji16}, the combination of several results from the literature implies that any such $\Gamma$ has property (T) and satisfies $\text{H}^2(\Gamma,\mathbb Z\Gamma)\not=0$.

Combined with a 2-cohomology version of a theorem by Moore and Schmidt \cite{MS80} which gives sufficient conditions for untwisting certain cocycles (see Section~\ref{ssec:MS}), Theorem~\ref{C} also leads to the following corollary.

\begin{mcor}\label{E}
Let $\Gamma$ be a countable group with property (T) that admits an ergodic p.m.p. action $\Gamma\curvearrowright (X,\mu)$ such that $\emph{H}^2(\Gamma,\emph{L}^0(X,\mathbb R))\not=0$.
In particular, assume that $\emph{H}^2(\Gamma,\mathbb R)\not=0$. 

Then C$^*(\Gamma)$ does not have the \emph{LP}. 
\end{mcor}

Corollary \ref{E} implies failure of the \text{LP} for several additional concrete classes of full group C$^*$-algebras.
Before discussing these classes, we record two well-known permanence properties for the (L)LP.

\begin{remark}\label{hereditary}
Let $\Sigma<\Gamma$ be countable groups. Then there are a canonical inclusion $C^*(\Sigma)\subset C^*(\Gamma)$ and conditional expectation $E:C^*(\Gamma)\rightarrow C^*(\Sigma)$ (see, e.g., \cite[Proposition 3.5]{Pi20}). This fact implies that if $C^*(\Gamma)$ has the LP (respectively, the LLP), then so does $C^*(\Sigma)$. Assume now that the inclusion $\Sigma<\Gamma$ has finite index. 
Then there are a $*$-isomorphism $\ell^{\infty}(\Gamma/\Sigma)\rtimes_{\text{f}}\Gamma\cong \mathbb B(\ell^2(\Gamma/\Sigma))\otimes C^*(\Sigma)$, a canonical inclusion $C^*(\Gamma)\subset\ell^{\infty}(\Gamma/\Sigma)\rtimes_{\text{f}}\Gamma$ and conditional expectation $E':\ell^{\infty}(\Gamma/\Sigma)\rtimes_{\text{f}}\Gamma\rightarrow C^*(\Gamma)$. From this it follows that if $C^*(\Sigma)$ has the LP (respectively, the LLP), then so does $C^*(\Gamma)$. 
\end{remark}

\begin{example}\label{example}(Full group C$^*$-algebras without the LP).

\begin{enumerate}[label=(\roman*)]
\item Let $G$ be a simple Lie group with trivial center, infinite cyclic fundamental group and property (T). This holds if $G=\text{Sp}_{2n}(\mathbb R)$, for $n\geq 2$, see \cite[Remark 3.5.5]{BdHV08} for more examples. Let $\Gamma<G$ be any lattice. For instance, take $\Gamma=\text{Sp}_{2n}(\mathbb Z)<G=\text{Sp}_{2n}(\mathbb R)$, for $n\geq 2$.  Then \cite[Corollary 3.5.6]{BdHV08} implies that $\text{H}^2(\Gamma,\mathbb Z)\not=0$. Using this, we observe in Lemma \ref{BdHV} that moreover $\text{H}^2(\Gamma,\mathbb R)\not=0$. Thus, $C^*(\Gamma)$ fails the \text{LP} by Corollary \ref{E}.
\item For a prime $p$, let $\Gamma_p$ be the kernel of the usual homomorphism $\text{SL}_3(\mathbb Z)\rightarrow\text{SL}_3(\mathbb Z/p\mathbb Z)$. By \cite{So78}, 
although $\text{H}^2(\text{SL}_3(\mathbb Z),\mathbb R)=0$, one has
$\text{H}^2(\Gamma_p,\mathbb R)\not=0$, for large $p$. Corollary \ref{E} implies that $C^*(\Gamma_p)$ does not have the \text{LP}, for every large enough prime $p$. In combination with Remark \ref{hereditary} this gives another proof of the failure of the \text{LP}  for $C^*(\text{SL}_3(\mathbb Z))$.
\item Let $\mathbb G$ be a simply connected, simple algebraic group over a number field $K$ with $K$-$\text{rank}$ at least $2$. Let $\mathcal O$ be the ring of integers of $K$ and $\Gamma =\mathbb G(\mathcal O)$ be the group of integral points in $\mathbb G$. Bekka \cite{Be99} proved that $C^*(\Gamma)$ is not residually finite dimensional. The structure theory of algebraic groups (see \cite{Be99}) implies that $\Gamma$ admits a subgroup that is isomorphic to a congruence subgroup of either $\text{SL}_3(\mathbb Z)$ or $\text{Sp}_4(\mathbb Z)$. Since $C^*(\text{SL}_3(\mathbb Z))$ and $C^*(\text{Sp}_4(\mathbb Z))$ do not have the \text{LP} by (i) and (ii), Remark \ref{hereditary} implies that $C^*(\Gamma)$ does not have the \text{LP}.
\item Let $\Gamma$ be a finitely presented residually finite group with property (T) and positive second $\ell^2$-Betti number, $\beta_2^{(2)}(\Gamma)>0$.  By L\"{u}ck's approximation theorem \cite{Lu94}, $\Gamma$ admits a finite index subgroup $\Sigma$ with positive second Betti number, $\beta_2(\Sigma)>0$, and thus $\text{H}^2(\Sigma,\mathbb R)\not=0$. Since $\Sigma$ has property (T), Corollary \ref{E} implies that $C^*(\Sigma)$ and $C^*(\Gamma)$ do not have the \text{LP}. \end{enumerate}
\end{example}

Next, motivated by Corollary \ref{E}, we introduce the following:

\begin{definition}
We denote by $\mathcal C$ the class of countable property (T) groups $\Gamma$ that admit an ergodic p.m.p. action $\Gamma\curvearrowright (X,\mu)$ such that $\text{H}^2(\Gamma,\text{L}^0(X,\mathbb R))\not=0$. 
\end{definition}

\begin{remark}\label{H^2(R)}
Using ergodic decomposition of p.m.p. actions, one can show that a property (T) group $\Gamma$ belongs to $\mathcal C$ if and only if it admits a p.m.p. (but not necessarily ergodic) action $\Gamma\curvearrowright (X,\mu)$ such that $\text{H}^2(\Gamma,\text{L}^0(X,\mathbb R))\not=0$ (see \cite[Proposition 2.5]{Ki17}). 
Note that $\mathcal C$ obviously contains any property (T) group $\Gamma$ with $\text{H}^2(\Gamma,\mathbb R)\not=0$, see Example \ref{example} for several families of such groups.
Additionally, it follows from \cite{Ji16} that any property (T) group $\Gamma$ with $\text{H}^2(\Gamma,\mathbb Z\Gamma)\not=0$ belongs to $\mathcal C$.
But besides these examples, not much is known about class $\mathcal C$. In fact, it  is unknown if $\mathcal C$ contains all property (T) groups.
Indeed, similar to the case of $\mathbb T$-valued second cohomology discussed in Remark \ref{H^2(T)}, no example of a free ergodic p.m.p. action $\Gamma\curvearrowright (X,\mu)$ of a property (T) group $\Gamma$ such that $\text{H}^2(\Gamma,\text{L}^0(X,\mathbb R))=0$ is known. 
\end{remark}

In Lemma \ref{embed}, we prove that the natural homomorphism $\text{H}^2(\Gamma,\text{L}^0(X,\mathbb R))\rightarrow\text{H}^2(\Gamma,\text{L}^0(X\times Y,\mathbb R))$ is injective, for any ergodic p.m.p. actions $\Gamma\curvearrowright (X,\mu)$ and $\Gamma\curvearrowright (Y,\nu)$ of a property (T) group $\Gamma$.
This implies that if $\Gamma\in\mathcal C$, then in fact there is a free ergodic p.m.p. action $\Gamma\curvearrowright (Z,\rho)$ with $\text{H}^2(\Gamma,\text{L}^0(Z,\mathbb R))\not=0$.
Moreover, it allows us to derive the following:

\begin{mprop}\label{ME}
The class $\mathcal C$ is closed under measure equivalence.
\end{mprop}
For the definition of Gromov's notion of measure equivalence, see \cite[Definition 1.1]{Fu99a}.
Recall that lattices in the same locally compact second countable group are measure equivalent (see \cite{Fu99a}). By combining Example \ref{example} (ii) with Proposition \ref{ME} and Corollary \ref{E}, it follows that $C^*(\Gamma)$ does not have the LP, for any lattice $\Gamma$ in $\text{SL}_3(\mathbb R)$. 

\begin{remark}
Let $\Gamma$ be a countable group and $\Gamma\curvearrowright (X,\mu)$ be a p.m.p. action. It would be interesting to determine if the natural homomorphism $\text{H}^2(\Gamma,\text{L}^2(X,\mathbb R))\rightarrow \text{H}^2(\Gamma,\text{L}^0(X,\mathbb R))$ is injective for property (T) groups. (As a positive indication that a result of this kind may hold, we show that the natural homomorphism $\text{H}^1(\Gamma,\text{L}^2(X,\mathbb R))\rightarrow \text{H}^1(\Gamma,\text{L}^0(X,\mathbb R))$ is injective, provided $\Gamma\curvearrowright (X,\mu)$ has spectral gap (see Theorem \ref{H1injective}).) 
If true, this would imply that $\mathcal C$ contains any property (T) group  $\Gamma$ that admits an orthogonal representation 
 $\pi:\Gamma\rightarrow\mathcal O(\mathcal H)$ on a real Hilbert space such that $\text{H}^2(\Gamma,\mathcal H)\not=0$. Indeed, for any orthogonal representation $\pi:\Gamma\rightarrow\mathcal O(\mathcal H)$, there is a p.m.p. action $\Gamma\curvearrowright (X_\pi,\mu_\pi)$ (called the Gaussian action associated to $\pi$) such that $\pi$ is contained in the Koopman representation of $\Gamma$ on $\text{L}^2(X,\mathbb R)$, and thus we have an embedding $\text{H}^2(\Gamma,\mathcal H)\subset\text{H}^2(\Gamma,\text{L}^2(X,\mathbb R))$.
 
\end{remark}

The  examples of property (T) groups without the LP provided above are typically finitely presented. On the other hand, we show that the LP fails for any infinitely presented property (T) group.

\begin{main}\label{infpres}
If $\Gamma$ is a non-finitely presented countable group with property (T), then $C^*(\Gamma)$ does not have the $\emph{LP}$.
\end{main}

Examples of non-finitely presented property (T) groups include $\text{SL}_3(\mathbb F_p[X])$, $\mathbb Z[1/p]^4\rtimes\text{Sp}_4(\mathbb Z[1/p])$, and infinite torsion quotients of uniform lattices in $\text{Sp}(n,1)$, 
where $p$ is a prime and $n\geq 2$, see \cite[Section 3.4]{BdHV08} and the references therein. 
Moreover, there are uncountably many pairwise non-isomorphic non-finitely presented property (T) groups, see \cite{Oz04a} and the references therein. 


In \cite[Corollary 5]{Oz04a}, Ozawa showed that within a certain uncountable family $\{\Gamma_{\alpha}\}_{\alpha\in \text{I}}$ of pairwise non-isomorphic property (T) groups, there is a group $\Gamma_\alpha$, for some $\alpha\in \text{I}$, such that $C^*(\Gamma_\alpha)$ fails the LP.
Theorem \ref{infpres} strengthens this result by showing that $C^*(\Gamma_\alpha)$ fails the LP, for every $\alpha\in \text{I}$.
 Indeed, by construction, the groups  $\{\Gamma_{\alpha}\}_{\alpha\in \text{I}}$ are all non-finitely presented. 

The proof of Theorem \ref{infpres} relies on a result of Shalom \cite{Sh00} asserting that every property (T) group is a quotient of a finitely presented property (T) group. The use of this result 
was inspired by an argument of Popa
 in \cite[Section 4]{Po07b}.

Next, we discuss a connection between our results and two recent notions of stability for groups. Given a tracial von Neumann algebra $(M,\tau)$, its $\text{L}^2$-norm is given by $\|x\|_{2,\tau}=\sqrt{\tau(x^*x)}$. If $M$ is a matrix algebra and $\tau$ its normalized trace, then $\|.\|_{2,\tau}$ is the normalized Hilbert-Schmidt norm.
Let $\Gamma$ be a countable group and $(M_n,\tau_n)$, $n\in\mathbb N$, be tracial von Neumann algebras. A sequence of maps $\varphi_n:\Gamma\rightarrow\mathcal U(M_n)$, $n\in\mathbb N$, is called an {\it asymptotic homomorphism} if it satisfies $\|\varphi_n(g)\varphi_n(h)-\varphi_n(gh)\|_{2,\tau_n}\rightarrow 0$, for every $g,h\in\Gamma$.

Following \cite[Definition~3]{HS18}, $\Gamma$ is called \textit{$W^*$-tracially stable} if for any sequence $(M_n,\tau_n)$, $n\in\N$, of tracial von Neumann algebras  and any asymptotic homomorphism $\varphi_n:\Gamma\rightarrow\cU(M_n)$, $n\in\mathbb N$, there exist homomorphisms $\psi_n:\Gamma\rightarrow \mathcal U(M_n)$,  $n\in\N$, such that $\|\varphi_n(g)-\psi_n(g)\|_{2,\tau_n}\rightarrow 0$, for every $g\in\Gamma$. If this holds whenever $M_n, n\in\mathbb N$, are matrix algebras, then $\Gamma$ is called {\it Hilbert-Schmidt stable} (or \textit{HS-stable}).

The proofs of our main results show that if $\Gamma$ is as in Theorem \ref{A} then $\Gamma$ is not HS-stable and   if $\Gamma$ is  as in Theorems \ref{C} or \ref{infpres} or  Corollaries \ref{D} or \ref{E}, then $\Gamma$ is not $W^*$-tracially stable. 
In fact, as explained below, we prove a much stronger statement: there is an asymptotic homomorphism $\varphi_n:\Gamma\rightarrow\mathcal U(M_n)$ (where $M_n$ are matrix algebras and tracial von Neumann algebras, respectively) for which no u.c.p. maps $\psi_n:C^*(\Gamma)\rightarrow M_n$ can be found such that $\|\varphi_n(g)-\psi_n(u_g)\|_2\rightarrow 0$, for every $g\in\Gamma$, where $\{u_g\mid g\in\Gamma\}$ are the canonical unitaries generating $C^*(\Gamma)$.
Thus, there are in fact no homomorphisms $\pi_n:\Gamma\rightarrow \mathcal U(M_n\bar{\otimes}\mathbb B(\ell^2))$ such that $\|\varphi_n(g)-(\pi_n(g))_{1,1}\|_{2,\tau_n}\rightarrow 0$, for every $g\in\Gamma$.
In particular,  any group $\Gamma$ satisfying the hypothesis of Theorem \ref{A} (e.g., $\Gamma=\mathbb Z^2\rtimes\text{SL}_2(\mathbb Z)$) is not flexibly HS-stable in the sense suggested in \cite[Section 4.4]{BL20}.
While it was shown in \cite{BL20} that infinite Connes-embeddable property (T) groups are not HS-stable, no examples of groups that are not flexibly HS-stable were previously known.


\subsection*{Outline of the proofs of Theorems \ref{A} and \ref{C}} We conclude the introduction with an outline of the proofs of our main technical results, Theorems \ref{A} and \ref{C}.
The starting point of our approach is the observation that if $C^*(\Gamma)$ has the LLP or the LP, then the set of asymptotic homomorphisms of $\Gamma$ is very rigid, see Section~\ref{ssec:consequence} (cf. \cite[Lemma 5]{Oz13}):

\begin{corollary}\label{ucplift}
Let $\Gamma$ be a countable group, $(M_n,\tau_n)$, $n\in\mathbb N$, be tracial von Neumann algebras and $\varphi_n:\Gamma\rightarrow\mathcal U(M_n)$, $n\in\mathbb N$, be an asymptotic homomorphism.

If  $C^*(\Gamma)$ has the \emph{LLP} and $M_n$ embeds into $R^{\omega}$, for every $n\in\mathbb N$, or $C^*(\Gamma)$ has the \emph{LP}, then there are u.c.p. maps $\psi_n:C^*(\Gamma)\rightarrow M_n$ such that $\|\varphi_n(g)-\psi_n(u_g)\|_{2,\tau_n}\rightarrow 0$, for every $g\in\Gamma$.
\end{corollary} 

Thus, to prove that $C^*(\Gamma)$ does not have the LP (respectively, the LLP), it suffices to produce an asymptotic homomorphism of $\Gamma$ to ($R^{\omega}$-embeddable) tracial von Neumann algebras which does not arise as an asymptotic perturbation of  a sequence of u.c.p. maps. However, in general, finding such asymptotic homomorphisms appears to be a non-trivial task. 

A main novelty of our approach is the use of asymptotic homomorphisms arising from 2-cocycles. We continue by explaining how this idea is implemented in the proofs of Theorems \ref{A} and \ref{C}.

First, in the context of Theorem \ref{A}, endow the matrix algebra $M_n:=\mathbb B(\mathcal H_n)$ with its normalized trace $\tau_n$. Then  
$\pi_n:\Gamma\rightarrow\mathcal U(M_n)$, $n\in\mathbb N$, is an asymptotic homomorphism. Assuming that $C^*(\Gamma)$ has the LLP, Corollary \ref{ucplift} gives u.c.p. maps $\psi_n:C^*(\Gamma)\rightarrow M_n$ so that $\|\pi_n(g)-\psi_n(u_g)\|_{2,\tau_n}\rightarrow 0$, for every $g\in\Gamma$. 

Next, we construct projective representations $\rho_n$ of $\Gamma$ whose $2$-cocycle is equal to $c_n$ and which have asymptotically invariant vectors. 
If the $\psi_n$'s would be $*$-homomorphisms,  we can define $\rho_n:\Gamma\rightarrow \mathcal U(\text{L}^2(M_n))$ by letting $\rho_n(g)(T)=\pi_n(g)T\psi_n(u_g)^*$ and note that the vectors $1\in\text{L}^2(M_n)$ are asymptotically invariant. 
In general, we define $\rho_n$ in a similar way via the Stinespring dilation of $\psi_n$.
On the other hand, \cite[Lemma 1.1]{NPS07} asserts that as the pair $(\Gamma,\Lambda)$ has relative property (T), given a projective representation $\rho$ of $\Gamma$ with almost invariant vectors, the restriction of its $2$-cocycle to $\Lambda$ must be a $2$-coboundary.
We thereby conclude that the restriction of $c_n$ to $\Lambda$ is a $2$-coboundary for large $n$. This gives a contradiction, implying that $C^*(\Gamma)$ does not have the LLP.

Second, we wish to adapt the proof of Theorem \ref{A} to prove Theorem \ref{C}.
In the setting of Theorem \ref{C}, let $M_n=\text{L}^{\infty}(X)\rtimes_{\sigma,c_n}\Gamma$ be the twisted crossed product von Neumann algebra associated to the action $\Gamma\curvearrowright^{\sigma}(X,\mu)$ and $2$-cocycle $c_n\in\text{Z}^2(\Gamma,\text{L}^0(X,\mathbb T))$. Let $\tau_n$ be the trace of $M_n$ and $\{u_{g,n}\mid g\in\Gamma\}\subset\mathcal U(M_n)$ be the canonical unitaries. Since $u_{g,n}u_{h,n}=c_n(g,h)u_{gh,n}$, for every $g,h\in\Gamma$, 
the maps $\Gamma\ni g\mapsto u_{g,n}\in\mathcal U(M_n)$, $n\in\mathbb N$, form an asymptotic homomorphism.  Assuming that $C^*(\Gamma)$ has the LP or that $C^*(\Gamma)$ has the LLP and the $M_n$'s are $R^{\omega}$-embeddable, Corollary \ref{ucplift} provides u.c.p. maps $\psi_n:C^*(\Gamma)\rightarrow M_n$ such that $\|u_{g,n}-\psi_n(u_g)\|_{2,\tau_n}\rightarrow 0$, for every $g\in\Gamma$. 

Suppose for the moment that the $\psi_n$'s are $*$-homomorphisms. Define $\rho_n:\Gamma\rightarrow\mathcal U(\text{L}^2(M_n))$ by letting  $\rho_n(g)(T)=u_{g,n}T\psi_n(u_g)^*$. Then $\rho_n$ satisfies $\rho_n(g)\rho_n(h)=c_n(g,h)\rho_n(gh)$, for every $g,h\in\Gamma$, and the vectors $1\in\text{L}^2(M_n)$ are asymptotically $\rho_n(\Gamma)$-invariant. In general, we first use a version of Stinespring's dilation theorem (see Section~\ref{ssec:dilation}) to dilate $\psi_n$ to a unital $*$-homomorphism from $C^*(\Gamma)$ to $M_n\bar{\otimes}\mathbb B(\ell^2)$, and then build the ``cocycle representations" $\rho_n$ in a similar fashion.

 Finally, we need to generalize the result from \cite{NPS07} to cocycle representations of $\Gamma$ whose cocycles can arise from a p.m.p. action of $\Gamma$. This is done in Section~\ref{ssec:NPSgen}. Altogether, these results allow us to once again deduce that the restriction of $c_n$ to $\Lambda$ has to be a $2$-coboundary for large $n$, leading to the desired contradiction.

\subsection*{Organization of the paper} Besides the introduction, this paper has seven other sections. In Section 2, we prove Corollary \ref{ucplift} and record various constructions and results concerning tracial von Neumann algebras, relative property (T) and group cohomology. Sections 3-8 are devoted to the proofs of our main results.

\subsection*{Acknowledgement} We would like to thank Sorin Popa for helpful comments.

\section{Preliminaries}
\subsection{Modules over tracial von Neumann algebras}\label{modules} 
In this paper, we will work with semifinite von Neumann algebras $\mathcal M$ endowed with a distinguished faithful normal semifinite trace $\text{Tr}$. We define the L$^2$-norm on $\mathcal M$ by $\|x\|_{2,\text{Tr}}=\sqrt{\text{Tr}(x^*x)}$, denote by L$^2(\mathcal M)$ the Hilbert space obtained by completing $\{x\in\mathcal M\mid \|x\|_{2,\text{Tr}}<\infty\}$ with respect to the L$^2$-norm and consider the standard representation $\mathcal M\subset\mathbb B($L$^2(\mathcal M))$.
In the particular case when $M$ is a finite von Neumann algebra with a faithful normal tracial state $\tau$, we call the pair $(M,\tau)$ a {\it tracial von Neumann algebra.}

Let $(M,\tau)$ be a tracial von Neumann algebra.
A Hilbert space $\mathcal H$ is  called a {\it left $M$-module} if it has a normal representation $\pi:M\rightarrow\mathbb B(\mathcal H)$ and a  {\it right M-module} if it has  a normal representation $\pi:M^{\text{op}}\rightarrow\mathbb B(\mathcal H)$. 
If $M$ is abelian, then the notions of left and right $M$-modules coincide, and we simply call them $M$-modules.
If $\mathcal H$ is a right $M$-module, then the conjugate Hilbert space $\overline{\mathcal H}=\{\overline{\xi}\mid\xi\in\mathcal H\}$ has a left $M$-module structure given by $x\overline{\xi}=\overline{\xi x^*}$.  

Let $\mathcal H$ be a right $M$-module.  
A vector $\xi\in\cH$ is called \textit{right bounded} if there is $C>0$ such that  $\norm{\xi x}\leq C\norm{x}_2$, for all $x\in M$. In other words, the map $M\ni x\mapsto \xi x\in\mathcal H$ extends to a bounded operator $L_\xi:$ L$^2(M)\rightarrow \cH$. The set of right bounded vectors $\cH_{\text{b}}$ forms a dense subspace of $\cH$. Moreover, for every $\xi,\eta\in \cH_{\text{b}}$, the operator $L_\eta^* L_\xi\in\mathbb B($L$^2(M))$ commutes with $M^{\text{op}}\subset\mathbb B($L$^2(M))$, and so belongs to $M$. We define an $M$-valued inner product on $\mathcal H_{\text{b}}$ by setting $(\xi,\eta)_M=L_\eta^* L_\xi$, and note that $\tau(x(\xi,\eta)_M)=\langle\xi x,\eta\rangle$, for every $x\in M$. 

If $\mathcal K$ is a left $M$-module, then the {\it Connes tensor product} $\mathcal H\otimes_M\mathcal K$ is the Hilbert space obtained by completing  $\mathcal H_{\text{b}}\otimes_{\text{alg}}\mathcal K/\ker (\langle\cdot,\cdot\rangle)$, where on the algebraic tensor product $\mathcal H_{\text{b}}\otimes_{\text{alg}}\mathcal K$ we consider the positive sesquilinear form given by $\langle\xi\otimes_M\eta,\xi'\otimes_M\eta'\rangle=\langle (\xi,\xi')_M\eta,\eta'\rangle$.

For further reference, we note that if $\xi,\eta\in\mathcal H_{\text{b}}$, then  in $\mathcal H\otimes_M\overline{\mathcal H}$ we have that \begin{equation}\label{estimate}\|\xi\otimes_M\overline{\eta}\|^2
=\langle(\xi,\xi)_M\eta,\eta\rangle=\tau((\xi,\xi)_M(\eta,\eta)_M)=\langle (\eta,\eta)_M\xi,\xi\rangle.\end{equation} %

Next, we describe a characterization of the Hilbert space $\cH\ot_M \bcH$ in terms of certain operators on $\cH$ that are ``Hilbert-Schmidt over $M$" (see \cite[Section~2]{Si11} and \cite[Chapter 8]{AP17}).
Let $\bB(\cH_M)\coloneqq (M^{\text{op}})'\cap\mathbb B(\mathcal H)$ be the algebra of right $M$-linear bounded operators on $\cH$. 
For $\xi,\eta\in\cH_b$ we define the ``rank one operator" $T_{\xi,\eta}=L_\xi L_\eta^*\in\mathbb B(\mathcal H_M)$.  

By \cite[Proposition 8.4.2]{AP17}, $\mathbb B(\mathcal H_M)$ admits a normal faithful semifinite trace $\widehat{\tau}$ such that  $$\widehat{\tau}(T_{\xi,\eta})=\tau((\xi,\eta)_M) = \langle \xi,\eta\rangle_\cH.$$
Then one can check that $\langle T_{\xi,\eta},T_{\xi',\eta'}\rangle_{\widehat{\tau}}=\langle\xi\otimes_M\overline{\eta},\xi'\otimes_M\overline{\eta'}\rangle$. Since the span of $\{T_{\xi,\eta}\mid \xi,\eta\in\mathcal H_{\text{b}}\}$ is a weak operator topology dense ideal of $\mathbb B(\mathcal H_M)$ (see \cite[Lemma 8.4.1]{AP17}), it follows that the map
\[
T_{\xi,\eta} \mapsto \xi\ot_M\bar\eta
\]
extends to a unitary isomorphism between L$^2(\mathbb B(\mathcal H_M),\widehat{\tau})$ and $\cH\ot_M \bcH$. 
When $M=\C$, this recovers the usual identification of $\cH\ot\bcH$ with the Hilbert space of Hilbert-Schmidt operators on $\cH$. 

\subsection{A consequence of the (L)LP for full group C$^*$-algebras}\label{ssec:consequence} 
In this section, we recall the definition of the LP and the LLP, and then prove Corollary \ref{ucplift}. Let $B$ be a C$^*$-algebra, $J\subset B$ a closed two-sided ideal, and $\pi:B\rightarrow B/J$ the quotient map, and $E$ an operator system. A u.c.p. map $\varphi:E\rightarrow B/J$ is {\it u.c.p. liftable} if there is a u.c.p. map $\psi:E\rightarrow B$ such that $\varphi=\pi\psi$.

\begin{definition}
	A C$^*$-algebra $A$ is said to have the {\it lifting property} (LP) if for every C$^*$-algebra $B$ and every closed two-sided ideal $J\subset B$, any u.c.p. map $\varphi:A\rightarrow B/J$ is u.c.p. liftable.
	A C$^*$-algebra $A$ has the {\it local lifting property} (LLP) if for every C$^*$-algebra $B$, every closed two-sided ideal $J\subset B$ and any u.c.p. map $\varphi:A\rightarrow B/J$, the restriction $\varphi_{|E}:E\rightarrow B/J$ of $\varphi$ to any finite dimensional operator system $E\subset A$ admits is u.c.p. liftable.
\end{definition}

A C$^*$-algebra $A$ has the {\it weak expectation property} (WEP) 
if there are a faithful representation $A\subset\mathbb B(H)$ and a u.c.p. map $\varphi:\mathbb B(H)\rightarrow A^{**}$ such that $\varphi_{|A}=\text{id}_A$. A C$^*$-algebra is called QWEP if it is a quotient of a C$^*$-algebra with the WEP. 
The proof of Corollary \ref{ucplift} relies on the following result due to Kirchberg \cite{Ki93}, see \cite[Corollary 3.12]{Oz04b}.
\begin{proposition}\label{lift}
	Let $A$ be a separable C$^*$-algebra with the \emph{LLP}, $B$ a \emph{QWEP} C$^*$-algebra and $J\subset B$ a closed two-sided ideal. Then any u.c.p. map $\varphi:A\rightarrow B/J$ is u.c.p. liftable.
\end{proposition}


\subsection*{Proof of Corollary \ref{ucplift}.} The proof is inspired by the proof of (2) $\Rightarrow$ (1) in \cite[Lemma 5]{Oz13}.

Consider the C$^*$-algebra $B=\prod_{n\in\mathbb N} M_n$, the closed ideal $J=\{(x_n)\in B\mid\lim\limits_{n\rightarrow\omega}\|x_n\|_{2,\tau_n}=0\}$ and the quotient map $\pi:B\rightarrow B/J$. Then $\varphi:\Gamma\rightarrow \mathcal U(B/J)$ given by $\varphi(g)=\pi((\varphi_n(g))_{n\in\mathbb N})$ is a homomorphism. Let $\Phi:C^*(\Gamma)\rightarrow B/J$  be the $*$-homomorphism given by $\Phi(u_g)=\varphi(g)$, for $g\in\Gamma$.

First, assume that $C^*(\Gamma)$ has the LP.  Then there is a u.c.p. map $\Psi:C^*(\Gamma)\rightarrow B$ such that $\Phi=\pi\Psi$. Hence we can find u.c.p. maps $\psi_n:C^*(\Gamma)\rightarrow M_n$, for every $n\in\mathbb N$, such that $\Psi=(\psi_n)_{n\in\mathbb N}$. Since $\pi((\varphi_n(g))_{n\in\mathbb N})=\Phi(u_g)=\pi((\psi_n(u_g))_{n\in\mathbb N})$, we get that $\lim\limits_{n\rightarrow\omega}\|\varphi_n(g)-\psi_n(u_g)\|_{2,\tau_n}=0$, for every $g\in\Gamma$. This finishes the proof in this case.

Secondly, assume that $C^*(\Gamma)$ has the {LLP} and $M_n$ embeds into $R^{\omega}$, for every $n\in\mathbb N$.
After replacing $M_n$ with the von Neumann subalgebra generated by $\{\varphi_n(g)\mid g\in\Gamma\}$, we may assume that $M_n$ is separable, for every $n\in\mathbb N$.
By a result of Kirchberg \cite{Ki93} (see \cite[Corollary 6.2]{Oz04b}), any separable tracial von Neumann algebra that embeds into $R^\omega$ is QWEP.  Since the class of QWEP C$^*$-algebras is closed under direct product (see \cite[Proposition 4.1]{Oz04b}),  we get that $B$ is QWEP. By Proposition \ref{lift}, there is a u.c.p. map $\Psi:C^*(\Gamma)\rightarrow B$ such that $\Phi=\pi\Psi$, and repeating the above argument finishes the proof.
\hfill$\blacksquare$

\subsection{Relative property (T)}\label{relT} Let $\Gamma$ be a countable group and $\Lambda$ be a subgroup.
The pair $(\Gamma,\Lambda)$ has {\it the relative property (T)} of Kazhdan-Margulis if there are a  finite set $F\subset\Gamma$ and $\delta>0$ such that if $\pi:\Gamma\rightarrow\mathcal U(\mathcal H)$ is a unitary representation and $\xi\in\mathcal H$ is a unit vector satisfying $\|\pi(g)\xi-\xi\|\leq\delta$, for every $g\in F$, then there is a non-zero vector $\eta\in \mathcal H$ that is $\pi(\Lambda)$-invariant, i.e., $\pi(h)\eta=\eta$, for every $h\in\Lambda$. 
The group $\Gamma$ has {\it property (T)} of Kazhdan if the pair $(\Gamma,\Gamma)$ has the relative property (T).

By a result of Jolissaint \cite{Jo05}, relative property (T) for $(\Gamma,\Lambda)$ is equivalent to the following formally stronger statement: for any $\eps>0$, there are a finite set $F=F(\varepsilon)\subset\Gamma$ and $\delta=\delta(\varepsilon)>0$ such that if $\pi:\Gamma\rightarrow \mathcal U(\cH)$ is a unitary representation  and $\xi\in\cH$ is a unit vector satisfying $\|\pi(g)\xi-\xi\|\leq\delta$, for every $g\in F$, then there is a $\pi(\Lambda)$-invariant vector $\eta\in \mathcal H$ such that $\|\eta-\xi\|\leq \eps$.

We next record a characterization of relative property (T) in terms of {projective} representations due to Nicoara, Popa and Sasyk \cite{NPS07}.
A  map $\pi:\Gamma\rightarrow\mathcal U(\mathcal H)$ is called a {\it projective representation} if $\pi(g)\pi(h)=c(g,h)\pi(gh)$, for some $c(g,h)\in\mathbb T$, for every  $g,h\in\Gamma$. 
Then $c:\Gamma\times\Gamma\rightarrow\mathbb T$ is a {\it 2-cocycle}, i.e., it satisfies $c(g,h)c(gh,k)=c(g,hk)c(h,k)$, for every $g,h,k\in\Gamma$.  A 2-cocycle $c$ is a {\it 2-coboundary} if there is a map $b:\Gamma\rightarrow\mathbb T$ such that $c(g,h)=b(g)b(h)b(gh)^{-1}$, for every $g,h\in\Gamma$. 

\begin{theorem}[{\!\!\!\cite[Lemma 1.1]{NPS07}}]\label{NPS}
If the pair $(\Gamma,\Lambda)$ has the relative property (T), then for any $\eps>0$, there are a finite set $F\subset\Gamma$ and $\delta>0$ such that the following holds.

Let  $\pi:\Gamma\rightarrow \mathcal U(\cH)$ be a projective representation with 2-cocycle $c:\Gamma\times\Gamma\rightarrow\mathbb T$  and $\xi\in\cH$ be a unit vector with $\inf\{\|\pi(g)\xi-\alpha\xi\|\mid \alpha\in\mathbb T\}<\delta$, for every $g\in F$. 
Then there are a map $b:\Lambda\rightarrow\mathbb T$ and $\eta\in \mathcal H$ such that  $c(h,k)=b(h)b(k)b(hk)^{-1}$ and $\pi(h)\eta=b(h)\eta$, for every $h,k\in\Lambda$, and $\|\eta-\xi\|<\eps$.
\end{theorem}

\subsection{Group cohomology} We next recall the notion of measurable group cohomology introduced by Moore in \cite{Mo76a}.
Let $G$ be a locally compact second countable group, $m_G$ a left invariant Haar measure of $G$, $A$ a Polish abelian group and  $\sigma:G\rightarrow\text{Aut}(A)$ be a continuous action of $G$ on $A$. 

For $n\geq 1$,  $\text{L}^0(G^n,A)$ denotes the set of equivalence classes of Borel maps $c:G^n\rightarrow A$, where two maps are equivalent if they are equal $m_G^n$-almost everywhere. Let $G^0=\{e\}$, so that $\text{L}^0(G^0,A)\equiv A$.
 Then $\text{L}^0(G^n,A)$ is a Polish group with respect to pointwise multiplication and the topology of convergence in measure (see \cite[page 6]{Mo76a}).
The {\it coboundary} map $\partial:\text{L}^0(G^n,A)\rightarrow \text{L}^{0}(G^{n+1},A)$  is defined by $$(\partial c)(g_1,...,g_{n+1})=\sigma_{g_1}(c(g_2,...,g_{n+1}))\cdot\prod_{i=1}^nc(g_1,...,g_ig_{i+1},...,g_{n+1})^{(-1)^i}\cdot c(g_1,...,g_n)^{(-1)^{n+1}}.$$

\begin{definition} Let $n\geq 1$.
Then $c\in\text{L}^0(G^n,A)$ is called an {\it n-cocycle} if $\partial c=e$ and an {\it n-coboundary} if there exists $b\in \text{L}^{0}(G^{n-1},A)$ such that $c=\partial b$.
Since $\partial(\partial c)=e$, any $n$-coboundary is an $n$-cocycle.

 We denote by $\text{Z}^n(G,A)$ and $\text{B}^n(G,A)$ the sets of $n$-cocycles and $n$-coboundaries, respectively. Then $\text{B}^n(G,A)\subset \text{Z}^n(G,A)$ are subgroups of $\text{L}^0(G^n,A)$, and the {\it n$^{th}$ measurable cohomology group of $G$ with coefficients in $A$} is defined as  $\displaystyle{\text{H}^n(G,A)=\frac{\text{Z}^n(G,A)}{\text{B}^n(G,A)}}.$
 Two $n$-cocycles $c,c'$ are said to be {\it cohomologous} if $cc'^{-1}$ is an $n$-coboundary. 
 \end{definition}

Note that $\text{Z}^n(G,A)$ is a closed subgroup of $\text{L}^0(G^n,A)$ and thus a Polish group. On the other hand, $\text{B}^n(G,A)$ is not necessarily closed.
 Nevertheless, as a consequence of basic facts about Borel functions and Polish groups (see, e.g., \cite[Theorem 1.4 and Theorem 2.6]{Ke94}), we have the following:
 

 \begin{lemma}\label{borel}
Let $n\geq 1$. Then $\emph{B}^n(G,A)$ is a Borel subgroup of $\emph{L}^0(G^n, A)$ and there is a Borel map
$\eta:\emph{B}^n(G,A)\rightarrow \emph{L}^{0}(G^{n-1},A)$ such that $\partial(\eta(c))=c$, for every $c\in \emph{B}^n(G, A)$.
\end{lemma}


\subsection{A vanishing result for symmetric 2-cocycles} We next record a result asserting that, for certain groups $A$, any 2-cocycle $c:\mathbb R\times \mathbb R\rightarrow A$  that is symmetric (i.e., satisfies $c(r,s)=c(s,r)$, for every $r,s\in\mathbb R$) must be a 2-coboundary.

\begin{lemma}\label{H^2} Let $A$ be either (1) a countable discrete abelian group, 
 or equal to (2) $\emph{L}^0(X,\mathbb T)$ or to (3) $\emph{L}^0(X,\mathbb T)/\mathbb T$, where $(X,\mu)$ is a probability space and we view $\mathbb T<\emph{L}^0(X,\mathbb T)$ as the closed subgroup of constant functions.
Then any symmetric 2-cocycle $c:\mathbb R\times\mathbb R\rightarrow A$ is a 2-coboundary.
\end{lemma}

{\it Proof.} Let $c:\mathbb R\times\mathbb R\rightarrow A$ be a symmetric 2-cocycle. By \cite[Theorem A]{AM13}, $c$ is cohomologous to a continuous 2-cocycle, so 
we may assume that $c$ is continuous. We discuss the three cases separately.

(1) Assume that $A$ is a countable discrete abelian group. Then since $\mathbb R$ is connected, $c$ must be a constant map, and thus it is a 2-coboundary.

(2)  Assume that $A=\text{L}^0(X,\mathbb T)$.  Define the group extension $E:=A\rtimes_{c}\mathbb R$ as the set $A\times\mathbb R$ with the group operation $(a,r)\cdot (b,s)=(c(r,s)ab,r+s)$. Then $E$ is a locally compact  second countable group and we have a short exact sequence $0\rightarrow A\rightarrow E\rightarrow \mathbb R\rightarrow 0$. Since $c$ is a symmetric $2$-cocycle, $E$ is abelian. By applying \cite[Theorem 4]{Mo76b} in the case $G=\mathbb R$, we conclude that $c$ is a $2$-coboundary.

 (3) Assume that $A=\text{L}^0(X,\mathbb T)/\mathbb T$. For $a\in \text{L}^0(X,\mathbb T)$, denote by $\bar{a}:=\mathbb T\cdot a$ its image in $A$. Let $i:A\rightarrow \text{L}^0(X\times X,\mathbb T)$ be given by $i(\bar a)(x,y)=a(x)a(y)^{-1}$. Since $i$ is a continuous homomorphism, $i\circ c:\mathbb R\times\mathbb R\rightarrow\text{L}^0(X\times X,\mathbb T)$ is a symmetric 2-cocycle. By applying (2), we get that $i\circ c$ is a 2-coboundary. Hence, there is a Borel function $b:\mathbb R\rightarrow\text{L}^0(X\times X,\mathbb T)$ such that $i(c_{r,s})=b_rb_sb_{r+s}^{-1}$, for $m_{\mathbb R}^2$-almost every $(r,s)\in\mathbb R^2$. Thus, we have
$$\text{$c_{r,s}(x)c_{r,s}(y)^{-1}=b_r(x,y)b_s(x,y)b_{r+s}(x,y)^{-1}$, for almost every $(r,s,x,y)\in \mathbb R^2\times X^2$.}$$ Hence, there is $y\in X$ such that $c_{r,s}(x)c_{r,s}(y)^{-1}=b_r(x,y)b_s(x,y)b_{r+s}(x,y)^{-1}$, for almost every $(r,s,x)\in\mathbb R\times X$. Thus, if we define $b:\mathbb R\rightarrow\text{L}^0(X,\mathbb T)$  by letting $b_r(x):=b_r(x,y)$, then we have that $c_{r,s}=\bar{b}_r\bar{b}_s\bar{b}_{r+s}^{-1}$, for almost every $(r,s)\in\mathbb R^2$. Hence, $c$ is a 2-coboundary.
\hfill$\blacksquare$

\subsection{Twisted group and crossed product von Neumann algebras}
Let $\Gamma$ be a countable group and $c\in\text{Z}^2(\Gamma,\mathbb T)$ be a $2$-cocycle. Then $\pi_c:\Gamma\rightarrow\mathcal U(\ell^2(\Gamma))$ given by $\pi_c(g)(\delta_h)=c(g,h)\delta_{gh}$ is a projective unitary representation with 2-cocycle $c$.
The {\it twisted group von Neumann algebra} $\text{L}_{c}(\Gamma)\subset\mathbb B(\ell^2(\Gamma))$ is defined as the weak operator closure of the linear span of $\pi_{c}(\Gamma)$. Then $\text{L}_c(\Gamma)$ is a tracial von Neumann algebra with its trace given by $\tau(\pi_c(g))=\delta_{g,e}$, for every $g\in\Gamma$.

More generally, let $\Gamma\curvearrowright^{\sigma} (A,\tau)$ be a trace preserving action of $\Gamma$ on an abelian tracial von Neumann algebra $(A,\tau)$ and $c\in\text{Z}^2(\Gamma,\mathcal U(A))$ be a $2$-cocycle. Define a map $\pi_c:\Gamma\rightarrow\mathcal U(\text{L}^2(A)\bar{\otimes}\ell^2(\Gamma))$ by letting 
$\pi_c(g)(x\otimes\delta_h)= \sigma_g(x)c(g,h)\otimes\delta_{gh}$, for every $g,h\in\Gamma$ and $x\in A$. We view $A$ as a  subalgebra of $\mathbb B(\text{L}^2(A)\bar{\otimes}\ell^2(\Gamma))$ by identifying it with $A\otimes 1$ and using the standard representation $A\subset\mathbb B(\text{L}^2(A))$. Then $\pi_c(g)\pi_c(h)=c(g,h)\pi_c(gh)$ and $\pi_c(g)x\pi_c(g)^*=\sigma_g(x)$, for all $g,h\in\Gamma$ and $x\in A$. The {\it twisted crossed product von Neumann algebra} $A\rtimes_{\sigma,c}\Gamma\subset\mathbb B(\text{L}^2(A)\bar{\otimes}\ell^2(\Gamma))$ is defined as the weak operator closure of the linear span of $\{a\pi_c(g)\mid a\in A, g\in\Gamma\}$. Then $A\rtimes_{\sigma,c}\Gamma$ is a tracial von Neumann algebra with its trace given by $\tau(a\pi_c(g))=\tau(a)\delta_{g,e}$, for every $a\in A$ and $g\in\Gamma$.


\section{Proofs of Theorem \ref{A} and Corollary \ref{B}}

In the first part of the section, we prove Theorem \ref{A} and Corollary \ref{B}. The proof of Theorem \ref{A} uses Corollary \ref{ucplift} to conclude that $C^*(\Gamma)$ does not have the LLP. In the second part of the section, we first give an alternative proof of this fact by showing that $C^*(\Gamma)\otimes_{\text{max}}\mathbb B(\ell^2)\not=C^*(\Gamma)\otimes_{\text{min}}\mathbb B(\ell^2)$, and then make a remark on the main result of \cite{JP95}. 

\subsection{ Proof of Theorem \ref{A}} 
Assume by contradiction that C$^*(\Gamma)$ has the LLP.
For every $n\in\mathbb N$, let
$M_n=\mathbb B(\mathcal H_n)$. We denote by  $\tau_n:M_n\rightarrow \mathbb C$ the normalized trace. 
Then for every $g,h\in\Gamma$, we have $\|\pi_n(g)\pi_n(h)-\pi_n(gh)\|_{2,\tau_n}=|c_n(g,h)-1|\rightarrow 0$.

Since $M_n$ is finite dimensional, Corollary \ref{ucplift} provides a u.c.p. map $\psi_n:C^*(\Gamma)\rightarrow M_n$, for every $n\in\mathbb N$, such that $\|\pi_n(g)-\psi_n(u_g)\|_{2,\tau_n}\rightarrow 0$, for every $g\in\Gamma$. By Stinespring's dilation theorem (see, e.g., \cite[Theorem 1.5.3]{BO08}), for every $n\in\mathbb N$ we can find a Hilbert space $\widehat{\mathcal H}_n\supset\mathcal H_n$ and a unitary representation $\varphi_n:\Gamma\rightarrow\mathcal U(\widehat{\mathcal H}_n)$ such that if $p_n:\widehat{\mathcal H}_n\rightarrow\mathcal H_n$ denotes the orthogonal projection, then $\psi_n(u_g)=p_n\varphi_n(g)p_n$, for all $g\in\Gamma$. Denote by $\text{Tr}:\mathbb B(\widehat{\mathcal H}_n)\rightarrow\mathbb C$ the usual semifinite trace.

Let $\mathcal K_n$ be the Hilbert space of bounded linear operators $T:\widehat{\mathcal H}_n\rightarrow \mathcal H_n$ endowed with the $\text{L}^2$-norm $\|T\|_{2,\text{Tr}}\coloneqq\sqrt{\text{Tr}(T^*T)}$.
Define $\rho_n:\Gamma\rightarrow\mathcal U(\mathcal K_n)$ by letting $\rho_n(g)(T)=\pi_n(g)T\varphi_n(g)^*$, for every $T\in \mathcal K_n$ and $g\in\Gamma$. Then $\rho_n$ is a projective representation with 2-cocycle equal to $c_n$.

Since $\displaystyle{\tau_n(x)=\frac{\text{Tr}(x)}{\text{Tr}(p_n)}}$, for all $x\in M_n$, for every $g\in\Gamma$ we get that 
\begin{align*}
\frac{\|\rho_n(g)(p_n)-p_n\|_{2,\text{Tr}}^2}{\|p_n\|_{2,\text{Tr}}^2}&= \frac{2\text{Tr}(p_n)-2\Re\;\text{Tr}(\pi_n(g)p_n\varphi_n(g)^*p_n)}{\text{Tr}(p_n)}\\&=2-2\Re\;\frac{\text{Tr}(\pi_n(g)\psi_n(u_g)^*)}{\text{Tr}(p_n)}\\&=2-2\Re\;\tau_n(\pi_n(g)\psi_n(u_g)^*)\\&=2\Re\;\tau_n(\pi_n(g)(\pi_n(g)-\psi_n(u_g))^*)\\&\leq 2\|\pi_n(g)-\psi_n(u_g)\|_{2,\tau_n}.
\end{align*}
Thus, the unit vectors $\xi_n=p_n/\|p_n\|_{2,\text{Tr}}\in \mathcal K_n$ are almost invariant: $\|\rho_n(g)(\xi_n)-\xi_n\|_{2,\text{Tr}}\rightarrow 0$, for every $g\in\Gamma$.
By applying Theorem \ref{NPS}, we conclude that for $n$ large enough, the restriction of $c_n$ to $\Lambda$ must be a 2-coboundary. This contradicts the hypothesis, and finishes the proof. \hfill$\blacksquare$

\subsection{Proof of Corollary \ref{B}} Let $\Gamma=\mathbb Z^2\rtimes\text{SL}_2(\mathbb Z)$ and $\Gamma_0=\mathbb Z^2\rtimes\Sigma$, where $\Sigma<\text{SL}_2(\mathbb Z)$ is a non-amenable subgroup.
Consider the 2-cocycle $c:\mathbb Z^2\times\mathbb Z^2\rightarrow \mathbb Z$ given by $c((x,y),(z,t))=xt-yz$. Then $c$ is $\text{SL}_2(\mathbb Z)$-invariant and thus it can be extended  to a 2-cocycle $c:\Gamma\times\Gamma\rightarrow\mathbb Z$ by letting $$\text{$c((a,g),(b,h))=c(a,g\cdot b)$, for all $a,b\in\mathbb Z^2$ and $g,h\in\text{SL}_2(\mathbb Z)$.}$$

For $k\in\mathbb N$, let $c_k\in\text{Z}^2(\Gamma,\mathbb T)$ be given by $c_k(g,h)=\exp(\frac{2\pi i\cdot c(g,h)}{k})$. First, we claim that the restriction of $c_k$ to $\mathbb Z^2$ is not a 2-coboundary, for any $k\geq 3$. Indeed, otherwise $c_k$ would be symmetric, i.e., $c_k(g,h)=c_k(h,g)$, for every $g,h\in\mathbb Z^2$. Letting $g=(1,0)$ and $h=(0,1)$ it would follow that $\exp(\frac{2\pi i}{k})=\exp(-\frac{2\pi i}{k})$ which would imply that $k\in\{1,2\}$. Second, note that $\lim\limits_{k\rightarrow\infty}c_k(g,h)=1$, for every $g,h\in\Gamma$. Third, let $\Gamma_k=(\mathbb Z/k\mathbb Z)^2\rtimes\text{SL}_2(\mathbb Z/k\mathbb Z)$ and $\theta_k:\Gamma\rightarrow\Gamma_k$ be the homomorphism given by reduction modulo $k$.
Then $c_k$ factors through $\theta_k\times\theta_k$, hence there is a 2-cocycle $\mu_k\in\text{H}^2(\Gamma_k,\mathbb T)$ such that $c_k=\mu_k\circ(\theta_k\times\theta_k)$. Let $\pi_k:\Gamma_k\rightarrow\mathcal U(\ell^2(\Gamma_k))$ be the projective representation associated to $\mu_k$. Then $\pi_k\circ\theta_k:\Gamma\rightarrow\mathcal U(\ell^2(\Gamma_k))$ is a projective representation with 2-cocycle $c_k$.

Since the pair $(\Gamma,\mathbb Z^2)$ has the relative property (T), we can apply Theorem \ref{A} to conclude that $C^*(\Gamma)$ does not have the LLP.
Moreover, since the pair $(\Gamma_0,\mathbb Z^2)$ has the relative property (T) by \cite{Bu91}, Theorem \ref{A} also implies that $C^*(\Gamma_0)$ does not have the LLP.
Since $\text{SL}_n(\mathbb Z)$ contains $\Gamma$ and $C^*(\Gamma)$ does not have the LLP,  Remark \ref{hereditary} implies that $C^*(\text{SL}_n(\mathbb Z))$ does not have the LLP, for any $n\geq 3$.

To prove the more general assertion, let $R$ be a finitely generated commutative ring with unit such that $\{2x\mid x\in R\}$ is infinite. 
Put $\Omega=R^2\rtimes\text{SL}_2(R)$. 
We may assume that $R$ has non-zero characteristic, for if $R$ has characteristic $0$, then $R$ contains $\mathbb Z$ as a subring and thus $\Omega$ contains $\mathbb Z^2\rtimes\text{SL}_2(\mathbb Z)$ as a subgroup. Then any finitely generated subgroup of $(R,+)$  is a torsion abelian group and thus finite. Moreover, since $R$ is a finitely generated commutative ring, a result of Baumslag \cite[Theorem 5.3]{Ba71} shows that $R$ is residually finite: for any finite set $F\subset R\setminus\{0\}$, there is an ideal $I$ of $R$ such that the quotient ring $R/I$ is finite and $F\cap I=\emptyset$. %

We will show that $\Omega$ satisfies the hypothesis of Theorem \ref{A} by adapting the above proof in the case $R=\mathbb Z$.
First, note that as above, the $2$-cocycle $c:R^2\times R^2\rightarrow R$ given by $c((x,y),(z,t))=xt-yz$ is $\text{SL}_2(R)$-invariant and thus it extends to a $2$-cocycle $c:\Omega\times\Omega\rightarrow R$. 

Let $\{a_n\mid n\in\mathbb N\}$ be an enumeration of $R$. Let $k\in\mathbb N$ and denote by $A_k$ the finite subgroup of $(R,+)$ generated by $\{a_1,...,a_k\}$. Since the set $\{2x\mid x\in R\}$ is infinite, we can find $x_k\in R$ such that $2x_k\notin A_k$. Let $B_k$ be the finite subgroup of $(R,+)$ generated by $A_k$ and $x_k$. 
 Since $R$ is residually finite, we can find an ideal $I_k$ of $R$ such that $R/I_k$ is finite and $B_k\cap I_k=\{0\}$. Let $p_k:R\rightarrow R/I_k$ be the quotient homomorphism. 
Since $B_k\cap I_k=\{0\}$, we get that $p_k(2x_k)\not\in p_k(A_k)$. Thus, we can find a character $\eta_k:R/I_k\rightarrow \mathbb T$ of the finite abelian group $(R/I_k,+)$ such that $\eta_k\equiv 1$ on $p_k(A_k)$ and $\eta_k(p_k(2x_k))\not=1$.
Then the character $\varphi_k:=\eta_k\circ p_k:(R,+)\rightarrow\mathbb T$ satisfies $\varphi_k\equiv 1$ on $A_k$ and $\varphi_k(2x_k)\not=1$.
Define $c_k\in\text{Z}^2(\Omega,\mathbb T)$ by letting $c_k=\varphi_k\circ c$. Since $c((x_k,0),(0,1))=x_k$ and $c((0,1),(x_k,0))=-x_k$, we get that $c_k((x_k,0),(0,1))=\varphi_k(x_k)\not=\varphi_k(-x_k)=c_k((0,1),(x_k,0))$. Thus, the restriction of $c_k$ to $R^2$ is not symmetric and hence not a $2$-coboundary, for every $k\in\mathbb N$.  Since $\varphi_k(a_1)=...=\varphi_k(a_k)=1$, we get that $\lim\limits_{k\rightarrow\infty}\varphi_k(a)=1$,  for every $a\in R$. Therefore, $\lim\limits_{k\rightarrow\infty}c_k(g,h)=1$, for every $g,h\in\Omega$. 
Moreover,  if we define $\Omega_k=(R/I_k)^2\rtimes\text{SL}_2(R/I_k)$ and denote by $\rho_k:\Omega\rightarrow\Omega_k$ the quotient homomorphism, then $c_k$ factors through $\rho_k\times\rho_k$.
Finally, by a result of Shalom \cite[Corollary 3.5]{Sh99}, the pair $(\Omega, R^2)$ has the relative property (T). Altogether, we can apply Theorem \ref{A} to conclude that $C^*(\Omega)$ does not have the LLP.

To prove the moreover assertion, let $\Delta$ be a countable group that contains $\Sigma=\mathbb F_2$ as a subgroup. 
Fix an embedding of $\Sigma$ into $\text{SL}_2(\mathbb Z)$ and let $A=\bigoplus_{\Delta/\Sigma}\mathbb Z^2$. Then the action of $\Sigma$ on $\mathbb Z^2$  can be used to build, via a construction called co-induction, an action of $\Delta$ on $A$ by automorphisms such that $A\rtimes\Delta$ contains $\mathbb Z^2\rtimes\Sigma$ as a subgroup (see the proof of \cite[Proposition 4.5]{Io11}). Since $C^*(\mathbb Z^2\rtimes\Sigma)$ does not have the LLP by the above, Remark \ref{hereditary} gives that $C^*(A\rtimes\Delta)$ does not have the LLP.
\hfill$\blacksquare$

\subsection{A second proof of Theorem \ref{A}} Next, we give an alternative proof of Theorem \ref{A} following closely \cite{Pi06}. 
Since the pair $(\Gamma,\Lambda)$ has the relative property (T), Theorem \ref{NPS} provides a finite set $F\subset\Gamma$ and $\delta>0$ such that  if $\pi:\Gamma\rightarrow \mathcal U(\mathcal H)$ is any projective representation whose $2$-cocycle $c\in \text{Z}^2(\Gamma,\mathbb T)$  satisfies $c_{|\Lambda}\not\in \text{B}^2(\Lambda,\mathbb T)$, then $\max_{g\in F}\|\pi(g)\xi-\xi\|\geq\delta\|\xi\|$, for any $\xi\in\mathcal H$. Thus, if we enumerate $F\cup\{e\}=\{g_1,...,g_m\}$, then
we have the following claim:
\begin{claim}\label{gap}  There is a constant $D<m$ such that
 if $\pi:\Gamma\rightarrow\mathcal U(\mathcal H)$ is any projective representation whose associated $2$-cocycle $c\in \emph{Z}^2(\Gamma,\mathbb T)$  satisfies $c_{|\Lambda}\not\in \emph{B}^2(\Lambda,\mathbb T)$ then
$\|\sum_{i=1}^m\pi(g_i)\|\leq D.$
\end{claim}

Indeed, for every $\xi\in\mathcal H$ we have
$$\delta^2\|\xi\|^2\leq\sum_{i,j=1}^m\|\pi(g_i)\xi-\pi(g_j)\xi\|^2=2m^2\|\xi\|^2-2\|\sum_{i=1}^m\pi(g_i)\xi\|^2.$$
This implies that $\|\sum_{i=1}^m\pi(g_i)\|\leq D$, where $D:=\sqrt{m^2-\frac{\delta^2}{2}}<m$, which proves the claim.

Let $\pi_n:\Gamma\rightarrow\mathcal U(\mathcal H_n)$, $n\in\mathbb N$, be projective representations as in the hypothesis of Theorem \ref{A}. 
Let $\mathbb B=\prod_{n\in\mathbb N}\mathbb B(\mathcal H_n)$ and $v_i=(\pi_n(g_i))_{n\in\mathbb N}\in \mathcal U(\mathbb B)$, for every $i\in\{1,...,m\}$. Define $$t=\sum_{i=1}^mu_{g_i}\otimes \overline{v}_i\in C^*(\Gamma)\otimes \overline{\mathbb B}.$$

Since the map $\Gamma\ni g\mapsto u_g\otimes\overline{\pi_n(g)}$ is a projective representation whose 2-cocycle $\overline{c_n}$ satisfies $\overline{c_n}_{|\Lambda}\notin\text{B}^2(\Lambda,\mathbb T)$, Claim \ref{gap} implies that $\|\sum_{i=1}^m u_{g_i}\otimes\overline{\pi_n(g_i)}\|\leq D$,
for every $n\in\mathbb N$. Thus, we get
\begin{equation}\label{min}
\|t\|_{C^*(\Gamma){\otimes}_{\text{min}}\overline{\mathbb B}}=\sup_{n\in\mathbb N}\|\sum_{i=1}^m u_{g_i}\otimes\overline{\pi_n(g_i)}\|\leq D.
\end{equation}

On the other hand,  we will show that $\|t\|_{C^*(\Gamma){\otimes}_{\text{max}}\overline{\mathbb B}}=m$ using the ultraproduct construction for tracial von Neumann algebras. 
To this end,  let $\tau_n:\mathbb B(\mathcal H_n)\rightarrow\mathbb C$ be the normalized trace and let $J=\{(x_n)_{n\in\mathbb N}\in\mathbb B\mid\lim\limits_{n\rightarrow\omega}\|x_n\|_{2,\tau_n}=0\}$. Then $J$ is a closed ideal of $\mathbb B$ and $M=\mathbb B/J$ is a tracial von Neumann algebra, with  its faithful normal tracial state $\tau$ given by $\tau((x_n)_{n\in\mathbb N})=\lim\limits_{n\rightarrow\omega}\tau_n(x_n)$.   Consider the $*$-homomorphisms $\lambda:M\rightarrow\mathbb B(\text{L}^2(M))$ and $\rho:\overline{M}\rightarrow\mathbb B(\text{L}^2(M))$ associated to left and right multiplication. 
Denote by $\eta:\mathbb B\rightarrow M$ the quotient map.

Since $\|\pi_n(g)\pi_n(h)-\pi_n(gh)\|_{2,\tau_n}\rightarrow 0$,
for every $g,h\in\Gamma$, we can define a unital $*$-homomorphism $\Phi:C^*(\Gamma)\rightarrow M$ by letting $\Phi(u_g)=\eta((\pi_n(g))_{n\in\mathbb N})$, for every $g\in\Gamma$. Since the $*$-homomorphisms $\lambda\circ\Phi:C^*(\Gamma)\rightarrow\mathbb B(\text{L}^2(M))$ and $\rho\circ\overline{\eta}:\overline{\mathbb B}\rightarrow\mathbb B(\text{L}^2(M))$ have commuting ranges we have that
\begin{equation}\label{max}
\|t\|_{C^*(\Gamma){\otimes}_{\text{max}}\overline{\mathbb B}}\geq \|(\lambda\circ\Phi\otimes \rho\circ\overline{\eta})(t)\|_{\mathbb B(\text{L}^2(M))}.
\end{equation}
Let $\xi\in \text{L}^2(M)$ be the unit vector corresponding to $1\in M$. Then $\Phi(u_{g_i})=\eta(v_i)$, $\overline{\eta}(\overline{v_i})=\overline{\eta(v_i)}$, thus $\lambda(\Phi(u_{g_i}))\rho(\overline{\eta}(\overline{v}_i))\xi=\eta(v_i)\xi\eta(v_i)^*=\xi$, for all $i\in\{0,1,...,m\}$. Therefore, $(\lambda\circ\Phi\otimes \rho\circ\overline{\eta})(t)\xi=m\xi$ and so  $\|(\lambda\circ\Phi\otimes \rho\circ\overline{\eta})(t)\|_{\mathbb B(\text{L}^2(M))}\geq m$. Together with \eqref{max} this gives that $\|t\|_{C^*(\Gamma){\otimes}_{\text{max}}\overline{\mathbb B}}\geq m$. 

In combination with \eqref{min}, we derive that $C^*(\Gamma)\otimes_{\text{max}}\overline{\mathbb B}\not=C^*(\Gamma)\otimes_{\text{min}}\overline{\mathbb B}$. Since $\mathbb B$ is isomorphic to $\overline{\mathbb B}$ and $\mathbb B$ can be embedded into $\mathbb B(\ell^2)$ such that it is the range of a conditional expectation, it follows that $C^*(\Gamma)\otimes_{\text{max}}\mathbb B(\ell^2)\not=C^*(\Gamma)\otimes_{\text{min}}\mathbb B(\ell^2)$. Therefore,  by \cite{Ki93} (see also \cite[Corollary 3.17]{Oz04b}) $C^*(\Gamma)$ does not have the \text{LLP}.
\hfill$\blacksquare$

\subsection{A remark on Junge and Pisier's theorem} In \cite{JP95}, Junge and Pisier established the fundamental result that $\mathbb B(\ell^2)\otimes_{\text{max}}\mathbb B(\ell^2)\not=\mathbb B(\ell^2)\otimes_{\text{min}}\mathbb B(\ell^2).$ The proof relies on the following crucial fact: for every $m\geq 3$, there exists a constant $C(m)<m$ and  an $m$-tuple of $N_k\times N_k$ unitary matrices $(u_1(k),...,u_m(k))\in \text{U}(N_k)^m$, for every $k\in\mathbb N$, such that $$\sup_{k\not=k'}\|\sum_{i=1}^mu_i(k)\otimes\overline{u_i(k')}\|\leq C(m).$$

There are by now several proofs of this fact, using property (T), expanders, quantum expanders or random matrices, for which we refer the reader to \cite[Chapters 18 and 19]{Pi20}. 
Our goal here is to show that this fact can also be deduced from Theorem \ref{NPS}, as follows.

 Let $\Gamma=\mathbb Z^2\rtimes\text{SL}_2(\mathbb Z)$ and $\Lambda=\mathbb Z^2$. Then $\Gamma$ is  2-generated. Indeed, since $S=\begin{pmatrix}0&-1\\ 1&0 \end{pmatrix}$ and $T=\begin{pmatrix} 1&1\\0&1\end{pmatrix}$ generate $\text{SL}_2(\mathbb Z)$, it follows that $g_1=aS$, where $a=\begin{pmatrix}1\\0 \end{pmatrix}\in\Lambda$, and  $g_2=T$ generate $\Gamma$. 

Let  $m\geq 3$ and put $g_3=...=g_m=e$. 
Since the pair $(\Gamma,\Lambda)$ has relative property (T) and  $\{g_1,g_2\}$ generates $\Gamma$, we get that there is a constant $D(m):=D<m$ such that Claim \ref{gap} holds.

Now, for $k\geq 2$, let $c_k\in\text{Z}^2(\Gamma,\mathbb T)$ and $\pi_k:\Gamma\rightarrow\mathcal U(\ell^2(\Gamma_k))$ be as defined in the proof of Corollary \ref{B}.
Then for every $k,k'\geq 2$ with $k\not= k'$, $\pi_k\otimes\overline{\pi_{k'}}:\Gamma\rightarrow\mathcal U(\ell^2(\Gamma_k)\otimes\overline{\ell^2(\Gamma_{k'})})$ is a projective representation  whose 2-cocycle $c_k\overline{c_{k'}}$ satisfies $(c_k\overline{c_{k'}})_{|\Lambda}\notin\text{B}^2(\Lambda,\mathbb T)$.
Thus by applying Claim \ref{gap} we get that $$\sup_{k\not=k'}\|\sum_{i=1}^m\pi_k(g_i)\otimes\overline{\pi_{k'}(g_i)}\|<D(m).$$

\section{Proof of Theorem \ref{C}}

This section is devoted to the proof of Theorem \ref{C}. 
\subsection{A 2-cohomology characterization of relative property (T)}\label{ssec:NPSgen} 
The first main ingredient in the proof of Theorem \ref{C} is a characterization of relative property (T) in terms of 2-cohomology. 

\begin{theorem}\label{NPSgen}
Let $\Gamma$ be a countable group and $\Lambda$ be a subgroup such that the pair $(\Gamma,\Lambda)$ has the relative property (T).  Then for every $\eps>0$, there exist $F\subset\Gamma$ finite and $\delta>0$ such that the following holds:

Let $\Gamma\curvearrowright^{\sigma} (A,\tau)$ be a trace preserving action of $\Gamma$ on an abelian tracial von Neumann algebra $(A,\tau)$ and $c\in\emph{Z}^2(\Gamma,\mathcal U(A))$ be a 2-cocycle. Let $\mathcal H$ be an $A$-module and $\pi:\Gamma\rightarrow \mathcal U(\mathcal H)$ be a map such that  $\pi(g)\pi(h)=c(g,h)\pi(gh)$ and $\pi(g)a\pi(g)^*=\sigma_g(a)$, for every $g,h\in\Gamma$ and $a\in A$.
Assume $\xi\in\mathcal H$ is a unit vector such that $(\xi,\xi)_A=1$ 
 and $\inf\{\|\pi(g)\xi-\xi a\|\mid a\in\mathcal U(A)\} <\delta$, for every $g\in F$.

Then there exist a $\sigma(\Lambda)$-invariant projection $p\in A$ with $\tau(p)> 1-\varepsilon$, a vector $\eta\in p\mathcal H$ with $(\eta,\eta)_A=p$ and a map $b:\Lambda\rightarrow\mathcal U(Ap)$ such that
\begin{enumerate}
\item $c(h,k)p=b_h\sigma_h(b_k)b_{hk}^*$, for every $h,k\in\Lambda$,
\item $\pi(h)(\eta)=b_h\eta$, for every $h\in\Lambda$, and
\item $\|\eta-\xi\|<\varepsilon$.
\end{enumerate} 
In particular, if $\eps<1$ and $\sigma_{|\Lambda}$ is ergodic, then the restriction of $c$ to $\Lambda$ must be a 2-coboundary.
\end{theorem}

In the case $A$ is equal to $\mathbb C1$ endowed with the trivial $\Gamma$-action, Theorem \ref{NPSgen} recovers Theorem \ref{NPS}.

{\it Proof.} After replacing $A$ with the smallest $\sigma(\Gamma)$-invariant von Neumann subalgebra that contains the image of $c$ and $\{a_g\mid g\in F\}$, where $a_g\in\mathcal U(A)$ is such that $\|\pi(a)\xi-\xi a_g\|<\delta$, we may assume that $A$ is separable.

Let $\eps>0$. Let $\eps_0>0$ such that $(14\eps_0)^2<\eps$ and $14\sqrt{2}\eps_0<\eps$. As $(\Gamma,\Lambda)$ has relative property (T), there are $F\subset\Gamma$ finite and $\delta>0$ such that if $\rho:\Gamma\rightarrow\mathcal U(\mathcal K)$ is a unitary representation and $\zeta\in\mathcal K$ is a unit vector satisfying $\|\rho(g)\zeta-\zeta\|<2\delta$, for every $g\in F$, then $\|\rho(h)\zeta-\zeta\|\leq\eps_0$, for every $h\in\Lambda$.

To prove that the conclusion holds for $F$ and $\delta$, consider the setting of Theorem \ref{NPSgen}.
Define the Connes tensor product $\cH \ot_A \bcH$. 
Since $(\pi(g)\zeta,\pi(g)\zeta')_A=\sigma_g((\zeta,\zeta')_A)$ and $\pi(g)\pi(h)=c(g,h)\pi(gh)$, for every $g,h\in\Gamma$ and $\zeta,\zeta'\in\mathcal H$, it follows that the map $\pi\ot\overline{\pi}:\Gamma\rightarrow \mathcal U(\cH\ot_A\bcH)$ given by $$(\pi\otimes\overline{\pi})(g)(\zeta\otimes_A\overline{\zeta'})=\pi(g)\zeta\otimes_A\overline{\pi(g)\zeta'}$$ is a well-defined unitary representation of $\Gamma$. 

Let $g\in F$. Then there is $a_g\in\mathcal U(A)$ such that $\|\pi(g)(\xi)-\xi a_g\|<\delta$.
Since $(\xi,\xi)_A=1$, we get that $(\pi(g)\xi,\pi(g)\xi)_A=(\xi a_g,\xi a_g)_A=1$ and \eqref{estimate} implies that
\begin{align*}\|(\pi\otimes\overline{\pi})(g)(\xi\otimes_A\overline{\xi})-\xi\otimes_A\overline{\xi}\|&=\|(\pi\otimes\overline{\pi})(g)(\xi\otimes_A\overline{\xi})-\xi a_g\otimes_A\overline{\xi a_g}\|
\\&\leq \|(\pi(g)(\xi)-\xi a_g)\otimes_A\overline{\pi(g)\xi}\|+\|\xi a_g\otimes_A\overline{(\pi(g)(\xi)-\xi a_g)}\|\\& = 2\|\pi(g)(\xi)-\xi a_g\|\\&<2\delta.\end{align*}

Next, by Section~\ref{modules}, there is a unitary $U:\cH\ot_A\bcH\rightarrow$ L$^2(\mathbb B(\mathcal H_A),\hat{\tau})$ such that $U(\xi\otimes_A\overline{\eta})=L_{\xi}L_{\eta}^*$, for every $\xi,\eta\in\mathcal H_{\text{b}}$.
Then we have that $U(\pi\ot\bar{\pi})(g)U^*(T) = \pi(g)T\pi(g)^*$, for every $g\in\Gamma$ and $T\in \mathbb B(\mathcal H_A)\cap \text{L}^2(\mathbb B(\mathcal H_A),\widehat{\tau}).$
Thus, if we let $Q=L_\xi L_\xi^*\in\mathbb B(\mathcal H_A)$ and use the above estimate, then $\|\pi(g)Q\pi(g)^*-Q\|_{2,\widehat{\tau}}<2\delta$, for all $g\in F$. 
By using relative property (T), we get that \begin{equation}\label{T}\text{$\|\pi(h)Q\pi(h)^*-Q\|_{2,\widehat{\tau}}\leq\eps_0$, for all $h\in\Lambda$. }\end{equation}

Since $L_{\xi}^*L_{\xi}=(\xi,\xi)_A=1$, we get that $Q$ is a projection. 
Let $\mathcal C\subset\mathbb B(\mathcal H_A)$ be the weak operator closure of the convex hull of the set $\{\pi(h)Q\pi(h)^*\mid h\in\Lambda\}$. 
By the proof of \cite[Lemma 14.3.3]{AP17}, $\mathcal C$ is $\|\cdot\|_{2,\widehat{\tau}}$-closed and the unique element $T$ of minimal $\|\cdot\|_{2,\widehat{\tau}}$-norm $T$ satisfies that \begin{equation}\label{TQ1}\text{$\|T\|_{\infty}\leq \|Q\|_{\infty}\leq 1$,\;\;\; $\|T-Q\|_{2,\widehat{\tau}}\leq\eps_0$\;\;\; and}\end{equation}
\begin{equation}\label{TQ2}\text{ $\pi(h)T\pi(h)^*=T$, for every $h\in\Lambda$.}
\end{equation}
By using that $Q$ is a projection and \eqref{TQ1} it follows that $\|TT^*-Q\|_{2,\widehat{\tau}}\leq 2\eps_0$, $\|(TT^*)^2-Q\|_{2,\widehat{\tau}}\leq 4\eps_0$ and $\|TT^*-(TT^*)^2\|_{2,\widehat{\tau}}\leq 6\eps_0$. Note that $|{\bf 1}_{[\frac{1}{2},\frac{3}{2}]}(t)-t|\leq 2|t-t^2|$, for every $t\geq 0$. Since $TT^*\in\mathbb B(\mathcal H_A)$ is a positive operator, if we define $P={\bf 1}_{[\frac{1}{2},\frac{3}{2}]}(TT^*)\in\mathbb B(\mathcal H_A)$ using  Borel functional calculus, then \begin{equation}\label{PQ}\|P-Q\|_{2,\widehat{\tau}}\leq \|P-TT^*\|_{2,\widehat{\tau}}+\|TT^*-Q\|_{2,\widehat{\tau}}\leq 2\|TT^*-(TT^*)^2\|_{2,\widehat{\tau}}+\|TT^*-Q\|_{2,\widehat{\tau}} \leq 14
\eps_0.\end{equation}
Moreover, \eqref{TQ2} implies that \begin{equation}\label{Pinv}\text{ $\pi(h)P\pi(h)^*=P$, for every $h\in\Lambda$.} \end{equation}
Since $(A,\tau)$ is a separable abelian tracial von Neumann algebra we can identify it with $($L$^{\infty}(X),\int\cdot\;\text{d}\mu)$, where $(X,\mu)$ is a standard probability space.  Since $\mathcal H$ is a right $A$-module, \cite[Theorem IV.8.21]{Ta01} implies that we can write $\mathcal H$ as a direct integral $\mathcal H=\int_X^{\oplus}\mathcal H_x\;\text{d}\mu(x)$ of Hilbert spaces $\{\mathcal H_x\mid x\in X\}$ such that  $a=\int_X^{\oplus}a(x)\cdot\text{Id}_{\mathcal H_x}\;\text{d}\mu(x)$, for every $a\in A$.

Let $\Gamma\curvearrowright (X,\mu)$ be the p.m.p. action such that $\sigma_g(a)(x)=a(g^{-1}x)$, for every $g\in\Gamma, a\in A$ and $x\in X$. Let $g\in\Gamma$. Since $\pi(g)a\pi(g)^*=\sigma_g(a)=a\circ g^{-1}$, for every $a\in A$, \cite[Theorem~IV.8.23]{Ta01} allows us to find a measurable field of unitary operators $\{\pi(g)_x:\mathcal H_x\rightarrow\mathcal H_{gx}\mid x\in X\}$, such that \begin{equation}\label{pi}\text{$\pi(g)=\int_X^{\oplus}\pi(g)_x\;\text{d}\mu(x)$, for every $g\in\Gamma$.}\end{equation} By applying \cite[Corollary IV.8.16]{Ta01} it follows that for every $S\in\mathbb B(\mathcal H_A)=A'\cap\mathbb B(\mathcal H)$, we can find an essentially bounded measurable field of operators $\{S_x:\mathcal H_x\rightarrow\mathcal H_x\mid x\in X\}$ such that $S=\int_X^{\oplus}S_x\;\text{d}\mu(x)$.  If we denote by $\text{Tr}$ the usual trace on $\mathbb B(\mathcal K)$, where $\mathcal K$ is a Hilbert space, then  \begin{equation}\label{trace}\text{$\widehat{\tau}(S)=\int_X\text{Tr}(S_x)\;\text{d}\mu(x)$, for every $S\in\mathbb B(\mathcal H_A), S\geq 0$.}\end{equation}

Let $X\ni x\mapsto \xi_x\in\mathcal H_x$ be a measurable field of vectors representing $\xi$. Since $(\xi,\xi)_A=1$ and $\|\xi_x\|^2=(\xi,\xi)_A(x)$, we may assume that $\xi_x$ is a unit vector for every $x\in X$.
Moreover, $Q=\int_X^{\oplus}Q_x\;\text{d}\mu(x)$, where $Q_x\in\mathbb B(\mathcal H_x)$ is the rank one orthogonal  projection onto $\mathbb C\xi_x$. Since $P\in\mathbb B(\mathcal H_A)$ is a projection, we can find a measurable field of projections $\{P_x\in\mathbb B(\mathcal H_x)\mid x\in X\}$ such that $P=\int_X^{\oplus}P_x\;\text{d}\mu(x)$. 

By combining \eqref{PQ} and \eqref{trace} we get that \begin{equation}\label{int} \int_X\|P_x-Q_x\|_{2,\text{Tr}}^2 \;\text{d}\mu(x)=\|P-Q\|_{2,\widehat{\tau}}^2\leq (14\varepsilon_0)^2,
\end{equation}
while the combination of \eqref{Pinv} and \eqref{pi} gives that \begin{equation}\label{invar} \text{$\pi(h)_xP_x\pi(h)_x^*=P_{hx}$, for every $h\in\Lambda$ and almost every $x\in X$.}
\end{equation}
Now, let $X_0=\{x\in X\mid\text{$P_x$ is a rank one orthogonal projection}\}$. 
Since for every $x\in X$ we have that $\|P_x-Q_x\|_{2,\text{Tr}}^2\geq |\text{Tr}(P_x)-\text{Tr}(Q_x)|=|\text{Tr}(P_x)-1|$, by using \eqref{int} we get that $\mu(X_0)\geq 1-(14\eps_0)^2$.
On the other hand, \eqref{invar} implies that $X_0$ is $\sigma(\Lambda)$-invariant. Let $X_0\ni x\mapsto \eta_x\in\mathcal H_x$ be a measurable field of unit vectors such that $P_x$ is the orthogonal projection onto $\mathbb C\eta_x$, for every $x\in X_0$. Moreover, we may assume that $\langle\xi_x,\eta_x\rangle\geq 0$, for every $x\in X_0$.

Let $h\in\Lambda$. Then \eqref{invar} implies that the orthogonal projections onto $\mathbb C\pi(h)_x(\eta_x)$ and $\mathbb C\eta_{hx}$ are equal for almost every $x\in X_0$. 
Thus, we can find a measurable map $b_h:X_0\rightarrow\mathbb T$ such that $\pi(h)_x(\eta_x)=b_h(hx)\eta_{hx}$, for almost every $x\in X_0$. 
 Let $\eta\in\int_{X_0}^{\oplus}\mathcal H_x\;\text{d}\mu(x)$ be the vector given by the measurable field $\{\eta_x\mid x\in X_0\}$. Identifying $\int_{X_0}^{\oplus}\mathcal H_x\;\text{d}\mu(x)$ with the $\pi(\Lambda)$-invariant subspace ${\bf 1}_{X_0}\mathcal H$ of $\mathcal H$, we get that $\pi(h)(\eta)=b_h\eta$. Since $\pi(h)\pi(k)=c(h,k)\pi(hk)$, we derive that $$
\text{$c(h,k)b_{hk}\eta=c(h,k)\pi(hk)(\eta)=\pi(h)\pi(k)\eta=\pi(h)(b_k\eta)=\sigma_h(b_k)b_h\eta$, for every $h,k\in\Lambda$.}$$

Finally, let $p={\bf 1}_{X_0}\in A$. Then $p\in A$ is a $\sigma(\Lambda)$-invariant projection with $\tau(p)\geq 1-(14\varepsilon_0)^2>1-\eps$,  $\eta\in p\mathcal H$ satisfies $(\eta,\eta)_A=p$ and condition (2), and we can view $b$ as a map $b:\Lambda\rightarrow\mathcal U(Ap)$.  Since $(\eta,\eta)_A=p$, the last displayed equation thus implies that $c(h,k)b_{hk}=\sigma_h(b_k)b_h$, for every $h,k\in\Lambda$, and condition (1) follows. To prove condition (3), note that if $x\in X_0$, then since $\langle \xi_x,\eta_x\rangle\geq 0$, we have $$\|P_x-Q_x\|_{2,\text{Tr}}^2=2(1-\langle \xi_x,\eta_x\rangle^2)\geq 2(1-\langle\xi_x,\eta_x\rangle)=\|\xi_x-\eta_x\|^2.$$
Thus, using \eqref{int}, we get that $\int_X\|\xi_x-\eta_x\|^2\;\text{d}\mu(x)=\int_{X_0}\|\xi_x-\eta_x\|^2\;\text{d}\mu(x)+\mu(X\setminus X_0)\leq 2(14\eps_0)^2$, which implies that $\|\xi-\eta\|\leq 14\sqrt{2}\eps_0<\eps$, as desired.
\hfill$\blacksquare$

\subsection{A dilation theorem}\label{ssec:dilation}

In the proof of Theorem \ref{C} we will also need a variant of Stinespring's dilation theorem.
 This result essentially goes back to \cite[Theorem 3]{Ka80} and is a one-sided version of a well-known dilation theorem for normal c.p. maps between tracial von Neumann algebras, see \cite[Chapter 1]{Po86}. For group C$^*$-algebras $B$ it appears as Theorem 2.3 in \cite{COT19}, where it is shown to hold for arbitrary von Neumann algebras $M$. 
 Nevertheless, we include the proof here for the reader's convenience.

\begin{proposition}\label{dilation}
Let $B$ be a unital C$^*$-algebra, $(M,\tau)$ be a tracial von Neumann algebra and $\varphi:B\rightarrow M$ be a u.c.p. map.

Then there are a Hilbert space $\mathcal H$, a projection $q\in \mathbb B(\mathcal H)\bar{\otimes}M$, a rank one projection $p\in\mathbb B(\mathcal H)$ and a $*$-homomorphism $\pi:B\rightarrow q(\mathbb B(\mathcal H)\bar{\otimes}M)q$ such that $p\otimes 1\leq q$ and $p\otimes\varphi(a)=(p\otimes 1)\pi(a)(p\otimes 1)$, for every $a\in B$. Moreover, we can take $q=1$.
\end{proposition}
{\it Proof.}
Consider the algebraic tensor product $B\ot_{\text{alg}} M$ equipped with the sesquilinear form
\[
\langle a\ot x,b\ot y\rangle = \tau(y^*\varphi(b^*a)x).
\]
Since $\varphi$ is a u.c.p. map, this form is positive definite and thus it gives an inner product on $(B\ot_{\text{alg}} M)/{\Ker(\langle .,.\rangle)}$. We denote by $\cK$ the completion of $(B\ot_{\text{alg}} M)/{\Ker(\langle .,.\rangle)}$.

We claim that there exist well-defined $*$-homomorphisms $\theta:B\rightarrow\mathbb B(\mathcal K)$ and $\rho:M^{\text{op}}\rightarrow\mathbb B(\mathcal K)$ 
such that for every  $a,b\in B$ and $x,y\in M$:
\[
\theta(a)(b\ot y)= ab\ot y \quad\text{and}\quad \rho(x)(b\ot y)= b\ot yx.
\]
The proof of this claim is standard (see, e.g., \cite[Section 13.1.2]{AP17}). Nevertheless, we recall the argument for completeness.
As $\varphi$ is u.c.p. we have $\sum_{i,j=1}^ny_j^*\varphi(b_j^*b_i)y_i\geq 0$, for all $b_1,...,b_n\in B$ and $y_1,...,y_n\in M$.
Thus, if $r=\sum_{i=1}^nb_i\otimes y_i\in B\otimes_{\text{alg}}M$, then for all $a\in B$ and $x\in M$ we have that
$$\text{$\langle\theta(a)r,\theta(a)r\rangle=\tau\Big(\sum_{i,j=1}^n y_j^*\varphi(b_j^*a^*ab_i)y_i\Big)\leq \|a^*a\|\cdot\tau\Big(\sum_{i,j=1}^n\tau(y_j^*\varphi(b_j^*b_i)y_i\Big)=\|a\|^2\cdot\langle r,r\rangle$}.$$
Also, by using that $\tau$ is a trace we derive that
$$\langle\rho(x)r,\rho(x)r\rangle=\tau\Big(xx^*\sum_{i,j=1}^ny_j^*\varphi(b_j^*b_i)y_i\Big)\leq \|x^*x\|\cdot \tau\Big(\sum_{i,j=1}^ny_j^*\varphi(b_j^*b_i)y_i\Big)=\|x\|^2\cdot\langle r,r\rangle.$$
These inequalities imply that $\theta$ and $\rho$ are indeed well-defined $*$-homomorphisms. 
Moreover, since the linear functional $M\ni x\mapsto \langle\rho(x)r,r\rangle=\tau(\big(\sum_{i,j=1}^ny_j^*\varphi(b_j^*b_i)y_i\big)x)$ is normal, it follows that $\rho$ is normal, or equivalently $\mathcal K$ is a right $M$-module.

Since $\varphi$ is unital, the map $V:\text{L}^2(M)\rightarrow\mathcal K$ given by $V(x)=1\otimes x$, for  $x\in M$, is a right $M$-modular isometry.
By the classification of right $M$-modules (see, e.g., \cite[Proposition 8.2.2]{AP17}), we can find a Hilbert space $\mathcal H$, a unit vector $\xi\in\mathcal H$ and a right $M$-modular isometry $W:\mathcal K\rightarrow\mathcal H\otimes \text{L}^2(M)$ such that $W(1\otimes x)=\xi\otimes x$, for all $x\in M$. Let $p\in\mathbb B(\mathcal H)$ be the rank one projection onto $\mathbb C\xi$ and $q\in\mathbb B(\mathcal H\otimes \text{L}^2(M))$ the orthogonal projection onto $W(\mathcal K)$. 
Then $p\otimes 1\leq q$. Moreover, since $W(\mathcal K)\subset\mathcal H\otimes \text{L}^2(M)$ is a right $M$-module, $q\in\mathbb B(\mathcal H)\bar{\otimes}M$.
We denote still by $W$ the unitary $W:\mathcal K\rightarrow W(\mathcal K)=q(\mathcal H\otimes \text{L}^2(M))$.

Let $\pi:B\rightarrow q\mathbb B(\mathcal H\otimes \text{L}^2(M))q$ be the $*$-homomorphism given by $\pi(a)=W\theta(a)W^*$. 
Since $\theta(B)$ commutes with $\rho(M)$ and $W$ is right $M$-modular,  $\pi(B)$ commutes with the right representation of $M$ on $q(\mathcal H\otimes L^2(M))$. Thus, $\pi(M)\subset q(\mathbb B(\mathcal H)\bar{\otimes}M)q$.

Finally, let $a\in B$. Since $W^*(\xi\otimes x)=1\otimes x$, for every $x,y\in M$ we have that
\begin{align*}\langle\pi(a)(\xi\otimes x),\xi\otimes y\rangle&=\langle \theta(a)W^*(\xi\otimes x),W^*(\xi\otimes y)\\&=\langle\theta(a)(1\otimes x),1\otimes y\rangle\\ &=\tau(y^*\varphi(a)x)\\&=\langle \xi\otimes\varphi(a)x,\xi\otimes y\rangle\end{align*} which implies that $(p\otimes 1)\pi(a)(p\otimes 1)=p\otimes\varphi(a)$. This proves the main assertion.

To justify the moreover assertion, note that after replacing $\mathcal H$ with a larger Hilbert space we may assume that the projections $1-q$ and $1$ are equivalent in $\mathbb B(\mathcal H)\bar{\otimes}M$ and there is a unital $*$-homomorphism $B\rightarrow \mathbb B(\mathcal H)$. Thus, $(1-q)(\mathbb B(\mathcal H)\bar{\otimes}M)(1-q)$ is $*$-isomorphic to $\mathbb B(\mathcal H)\bar{\otimes}M$ and hence there is a unital $*$-homomorphism $\pi':B\rightarrow (1-q)(\mathbb B(\mathcal H)\bar{\otimes}M)(1-q)$. 
Then the unital $*$-homomorphism $\pi\oplus\pi':B\rightarrow \mathbb B(\mathcal H)\bar{\otimes}M$ shows that we may take $q=1$.
\hfill$\blacksquare$

\subsection{Proof of Theorem \ref{C}} We prove the following stronger version of Theorem \ref{C}.

\begin{theorem}
Let $\Gamma$ be a countable group and  $\Lambda$ a subgroup such that $(\Gamma,\Lambda)$ has relative property (T).  Assume that there is a trace preserving action $\Gamma\curvearrowright^{\sigma} (A,\tau)$ on an abelian tracial von Neumann algebra $(A,\tau)$ such that $\sigma_{|\Lambda}$ is ergodic, and  2-cocycles $c_n\in\emph{Z}^2(\Gamma,\mathcal U(A))$ such that the restriction of $c_n$ to $\Lambda$ is not a 2-coboundary, for every $n\in\mathbb N$, and $\lim\limits_{n\rightarrow\infty}\|c_n(g,h)-1\|_{2,\tau}=0$, for every $g,h\in\Gamma$.

Then C$^*(\Gamma)$ does not have the \emph{LP}. 
Moreover, if  the twisted crossed product von Neumann algebra $A\rtimes_{\sigma,c_n}\Gamma$ embeds into $R^{\omega}$, for every $n\in\mathbb N$,
 then C$^*(\Gamma)$ does not have the \emph{LLP}.

\end{theorem}

{\it Proof.}
 For $n\in\mathbb N$, denote $M_n=A\rtimes_{\sigma,c_n}\Gamma$. Let $\{u_{g,n}\mid g\in\Gamma\}\subset\mathcal U(M_n)$ and $\tau_n:M_n\rightarrow\mathbb C$ be the canonical unitaries and trace. Then $$\text{$\|u_{g,n}u_{h,n}-u_{gh,n}\|_{2,\tau_n}=\|c_n(g,h)-1\|_{2,\tau}\rightarrow 0$, for every $g,h\in\Gamma$.}$$
Assuming by contradiction that the conclusion fails, Corollary \ref{ucplift} implies the existence of u.c.p. maps $\psi_n:C^*(\Gamma)\rightarrow M_n$, $n\in\mathbb N$, such that $\|u_{g,n}-\psi_n(u_g)\|_{2,\tau_n}\rightarrow 0$, for every $g\in\Gamma$.
By Proposition \ref{dilation}, we can find a Hilbert space $\mathcal H_n$, a rank one projection $p_n\in \mathbb B(\mathcal H_n)$ and a $*$-homomorphism $\rho_n:C^*(\Gamma)\rightarrow \mathbb B(\mathcal H_n)\bar{\otimes}M_n$ such that $p_n\otimes\psi_n(x)=(p_n\otimes 1)\rho_n(x)(p_n\otimes 1)$, for every $x\in C^*(\Gamma)$, $n\in\mathbb N$.

For $n\in\mathbb N$, consider the von Neumann algebra $\mathcal M_n=\mathbb B(\mathcal H_n)\bar{\otimes}M_n$ together with the semifinite trace $\widehat{\tau}_n=\text{Tr}\otimes\tau_n$, where $\text{Tr}$ denotes the usual trace on $\mathbb B(\mathcal H_n)$, and the associated L$^2$-norm, $\|\cdot\|_{2,\widehat{\tau}_n}$. For $x\in \mathcal M_n$ with $\|x\|_{2,\widehat{\tau}_n}<\infty$, we denote still by $x$ the corresponding vector in L$^2(\mathcal M_n)$. We  consider the standard unital normal embeddings $\mathcal M_n,\mathcal M_n^{\text{op}}\subset\mathbb B(\text{L}^2(\mathcal M_n))$ given by the left and right multiplication actions of $\mathcal M_n$. 
Using that $p_n\otimes 1\in\mathcal M_n$, we define $\mathcal K_n=(p_n\otimes 1)\text{L}^2(\mathcal M_n)$. 

Then we have unital  normal embeddings $p_n\otimes M_n=(p_n\otimes 1)\mathcal M_n(p_n\otimes 1)\subset\mathbb B(\mathcal K_n)$ and $\mathcal M_n^{\text{op}}\subset\mathbb B(\mathcal K_n)$.
obtained by noticing that the left multiplication action of $(p_n\otimes 1)\mathcal M_n(p_n\otimes 1)$  on L$^2(\mathcal M_n)$ and the right multiplication action of $\mathcal M_n$ on L$^2(\mathcal M_n)$ leave $\mathcal K_n$ invariant. 

In particular, we can endow $\mathcal K_n$ with a left $A$-module structure given by $a\cdot\xi=(p_n\otimes a)\xi$.
 Further, for $n\in\mathbb N$, we define a map $\pi_n:\Gamma\rightarrow \mathcal U(\mathcal K_n)$ by letting $$\text{$\pi_n(g)\xi=(p_n\otimes u_{g,n})\xi\rho_n(u_g)^*$, for every $g\in\Gamma$ and $\xi\in\mathcal K_n$.}$$
Using that $\rho_n$ is a homomorphism, $u_{g,n}u_{h,n}=c_n(g,h)u_{gh,n}$ and $u_{g,n}au_{g,n}^*=\sigma_g(a)$ we get that $\pi_n(g)\pi_n(h)=c_n(g,h)\pi_n(gh)$ and $\pi_n(g)a\pi_n(g)^*=\sigma_g(a)$, for every $g,h\in\Gamma$ and $a\in A$.

Next, let $\xi_n=p_n\otimes 1\in\mathcal K_n$. Then $\|\xi_n\|_{2,\widehat{\tau}_n}=1$. Since $\langle a\cdot\xi_n,\xi_n\rangle=\langle p_n\otimes a,p_n\otimes 1\rangle=\tau(a)$, for every $a\in A$, we get that $(\xi_n,\xi_n)_A=1$. Moreover, if $g\in\Gamma$, then we have \begin{align*}\langle\pi_n(g)\xi_n,\xi_n\rangle&=\langle(p_n\otimes u_{g,n})(p_n\otimes 1)\rho_n(u_g)^*,p_n\otimes 1\rangle\\&=\widehat{\tau}_n((p_n\otimes u_{g,n})(p_n\otimes 1)\rho_n(u_g)^*(p_n\otimes 1))\\&=\widehat{\tau}_n((p_n\otimes u_{g,n})(p_n\otimes\psi_n(u_g)^*))\\&=\tau_n(u_{g,n}\psi_n(u_g)^*).\end{align*} 
Thus, we get that \begin{align*}\|\pi_n(g)\xi_n-\xi_n\|_{2,\widehat{\tau}_n}^2&=2(1-\Re\langle\pi_n(g)\xi_n,\xi_n\rangle)\\&=2(1-\Re\tau_n(u_{g,n}\psi_n(u_g)^*))\\&=2\Re\tau_n(u_{g,n}(u_{g,n}-\psi_n(u_g))^*) \\&\leq 2\|u_{g,n}-\psi_n(u_g)\|_{2,\tau_n}.\end{align*}
Hence, $\|\pi_n(g)\xi_n-\xi_n\|_{2,\widehat{\tau}_n}\rightarrow 0$, for every $g\in\Gamma$.
Altogether, we are in position to apply Theorem \ref{NPSgen} and conclude that the restriction of $c_n$ to $\Lambda$ is a 2-coboundary, for large enough $n$. This contradicts the hypothesis and finishes the proof. \hfill$\blacksquare$

\section{Proof of Corollary \ref{D}}
Let $\Gamma$ be a property (T) group with $\text{H}^2(\Gamma,\mathbb Z\Gamma)\not=\{0\}$. Assume by contradiction that $C^*(\Gamma)$ has the \text{LP}.
Let $d\in \text{Z}^2(\Gamma,\mathbb Z\Gamma)\setminus\text{B}^2(\Gamma,\mathbb Z\Gamma)$. Define $\widetilde d=(d,-d):\Gamma\times\Gamma\rightarrow\mathbb Z\Gamma\oplus\mathbb Z\Gamma$. Then $\widetilde d\in\text{Z}^2(\Gamma,\mathbb Z\Gamma\oplus\mathbb Z\Gamma)$.


For a countable abelian group $G$, consider the embedding $i:G\rightarrow \text{L}^0(\widehat{G},\mathbb T)$ given by $i(a)(\varphi)=\varphi(a)$.
If $\Gamma$ acts on $G$ by automorphisms and we endow $\widehat{G}$ with the dual $\Gamma$-action, then $i$ is $\Gamma$-equivariant.

Next, 
let $\lambda$ be the Haar measure of $\mathbb T$ and consider the Bernoulli action $\Gamma\curvearrowright^{\sigma} (\mathbb T^{\Gamma},\lambda^{\Gamma})$. Denote by $\Gamma\curvearrowright^{\widetilde{\sigma}}(\mathbb T^{\Gamma}\times\mathbb T^{\Gamma},\lambda^{\Gamma}\times\lambda^{\Gamma})$ the diagonal action $\sigma^2$.
Note that the dual actions of $\Gamma$ on $\widehat{\mathbb Z\Gamma}=\mathbb T^{\Gamma}$ and $\widehat{\mathbb Z\Gamma\oplus\mathbb Z\Gamma}=\mathbb T^{\Gamma}\times\mathbb T^{\Gamma}$, are $\sigma$ and $\widetilde\sigma$, respectively.
As $i$ is $\Gamma$-equivariant,  following \cite{Ji16}, we define $$\text{$c:=i\circ d\in\text{Z}^2(\Gamma,\text{L}^0(\mathbb T^{\Gamma},\mathbb T))$\;\; and\;\; $\widetilde c:=i\circ\widetilde{d}\in\text{Z}^2(\Gamma,\text{L}^0(\mathbb T^{\Gamma}\times\mathbb T^{\Gamma},\mathbb T))$.}$$

Let $A=\text{L}^{\infty}(\mathbb T^{\Gamma})$ and denote also by $\sigma$ and $\widetilde{\sigma}$ the Bernoulli actions of $\Gamma$ on $A$ and $A\bar{\otimes}A$. With this notation, we can view $c\in\text{Z}^2(\Gamma,\mathcal U(A))$ and $\widetilde{c}\in\text{Z}^2(\Gamma,\mathcal U(A\bar{\otimes}A))$.
Since $\sigma$ is malleable in the sense of Popa (see, e.g., \cite{Po07a}), there is a 1-parameter group $\{\alpha_t\}_{t\in\mathbb R}$ of automorphisms of $A\bar{\otimes}A$ such that 
\begin{enumerate}
\item $\alpha_0=\text{Id}$, $\alpha_1(a\otimes 1)=1\otimes a$, for every $a\in A$, 
\item $\alpha_t$ commutes with $\widetilde{\sigma}$, for every $t\in\mathbb R$, and
\item $\lim\limits_{t\rightarrow 0}\|\alpha_t(a)-a\|_2=0$, for every $a\in A\bar{\otimes}A$.
\end{enumerate}
In particular, using (2), we can define $c_t:=c\;\!\alpha_t(c)^*\in\text{Z}^2(\Gamma, \mathcal U(A\bar{\otimes}A))$, for every $t\in\mathbb R$, where we identify $A$ with $A\otimes 1\subset A\otb A$. Then $\lim\limits_{t\rightarrow 0}\|c_t(g,h)-1\|_2=\lim\limits_{t\rightarrow 0}\|\alpha_t(c(g,h))-c(g,h)\|_2=0$ by (3), for every $g,h\in\Gamma$.  Since $C^*(\Gamma)$ is assumed to have the \text{LP}, Theorem \ref{C} implies that there is $\delta>0$ such that $c_t\in \text{B}^2(\Gamma,\mathcal U(A\bar{\otimes}A))$, for every $t\in\mathbb R$ with $|t|<\delta$. Thus, there is $n\in\mathbb N$ such that $c_{\frac{1}{n}}\in \text{B}^2(\Gamma,\mathcal U(A\bar{\otimes}A))$. Therefore \begin{equation}\label{cob}\widetilde c=c\otimes c^*=c\;\alpha_1(c)^*=c_{\frac{1}{n}}\alpha_{\frac{1}{n}}(c_{\frac{1}{n}})...\alpha_{\frac{n-1}{n}}(c_{\frac{1}{n}})\in\text{B}^2(\Gamma,\mathcal U(A\bar{\otimes}A)).\end{equation}
Since $\widetilde{\sigma}$ is a Bernoulli action, so are the diagonal actions ${\widetilde{\sigma}}^2$ and $\widetilde{\sigma}^4$. Since $\Gamma$ has property (T), Popa's cocycle superrigidity theorem \cite{Po07a} implies that ${\widetilde{\sigma}}^2$ and $\widetilde{\sigma}^4$ are $\{\mathbb T\}$-cocycle superrigid.
Since $\Gamma$ has property (T), \cite[Corollary 4.1]{Ji16} implies that the group $\text{H}^2(\Gamma,\mathbb Z\Gamma\oplus\mathbb Z\Gamma)\cong\text{H}^2(\Gamma,\mathbb Z\Gamma)\oplus \text{H}^2(\Gamma,\mathbb Z\Gamma)$ is torsion free. Finally, since $\widetilde{c}=i\circ\widetilde{d}$, \cite[Theorem 1.1]{Ji16} implies that $\widetilde{d}=(d,-d)\in\text{B}^2(\Gamma,\mathbb Z\Gamma\oplus\mathbb Z\Gamma)$. This gives that $d\in\text{B}^2(\Gamma,\mathbb Z\Gamma)$, which is a contradiction. 
\hfill$\blacksquare$

\section{Proof of Corollary \ref{E}}
In  this section we prove Corollary \ref{E} and then justify a claim made in Example \ref{example} (i).

\subsection{A Moore-Schmidt theorem for 2-cohomology}\label{ssec:MS} Let $\Gamma\curvearrowright (X,\mu)$ be a p.m.p. action and $c\in \text{Z}^1(\Gamma, \text{L}^0(X,\mathbb R))$ be a 1-cocycle. For $r\in\mathbb R$, define $c_r\in\text{Z}^1(\Gamma,\text{L}^0(X,\mathbb T))$ by $c_r(g)=\exp(irc(g))$. 
If $c$ is a 1-coboundary, then so is $c_r$, for every $r\in\mathbb R$.

Conversely, a theorem of Moore and Schmidt (see \cite[Theorem 4.3]{MS80}) shows that if $c_r$ is a 1-coboundary, for every $r\in\mathbb R$, then so is $c$.
The proof of Corollary \ref{E} relies on the following  analogue of this result  for 2-cohomology.

\begin{theorem}\label{untwist}
Let $\Gamma$ be a countable group and $\Gamma\curvearrowright (X,\mu)$ be an ergodic p.m.p. action. Let $c\in\emph{Z}^2(\Gamma, \emph{L}^0(X,\mathbb R))$ and for $r\in\mathbb R$, define $c_r\in\emph{Z}^2(\Gamma, \emph{L}^0(X,\mathbb T))$  by $c_r(g,h)(x)=\exp(irc(g,h)(x))$, for $g,h\in\Gamma$ and $x\in X$. 
Assume that ${c_r}\in \emph{B}^2(\Gamma,\emph{L}^0(X,\mathbb T))$, for every $r\in\mathbb R$. 

Suppose further that $\emph{B}^1(\Gamma, \emph{L}^0(X,\mathbb T))$ is an open subgroup of $\emph{Z}^1(\Gamma,\emph{L}^0(X,\mathbb T))$.

Then $c\in\emph{B}^2(\Gamma, \emph{L}^0(X,\mathbb R))$.
\end{theorem}

The condition that $\text{B}^1(\Gamma, \text{L}^0(X,\mathbb T))$ is an open subgroup of $\text{Z}^1(\Gamma,\text{L}^0(X,\mathbb T))$ is satisfied whenever $\Gamma$ has property (T) by a result of Schmidt
\cite[Theorem 3.4]{Sc81}. It also holds if the action $\Gamma\curvearrowright (X,\mu)$ is $\{\mathbb T\}$-cocycle superrigid, see Popa's cocycle superrigidity theorems \cite{Po07a,Po08} for such actions, and $\text{Char}(\Gamma)=\{e\}$.
Finally, the condition trivially holds if $X$ consists of one point, or equivalently when $\text{L}^0(X,\mathbb T)=\mathbb T$, leading to the following consequence:

\begin{corollary}\label{Rvalued}
Let $\Gamma$ be a countable group and $c\in \emph{Z}^2(\Gamma,\mathbb R)$. For $r\in\mathbb R$, define $c_r\in\emph{Z}^2(\Gamma,\mathbb T)$ by letting $c_r=\exp(irc)$.
Assume that ${c_r}\in\emph{B}^2(\Gamma,\mathbb T)$, for every $r\in\mathbb R$. Then $c\in\emph{B}^2(\Gamma,\mathbb R)$.
\end{corollary}

{\it Proof of Theorem~\ref{untwist}.} Put $A=\text{L}^0(X,\mathbb T)$.
Since the map $\mathbb R\ni r\mapsto c_r\in \text{Z}^2(\Gamma,A)$ is continuous and its image is contained in $\text{B}^2(\Gamma,A)$, Lemma \ref{borel} implies the existence of a Borel map $\mathbb R\ni r\mapsto b_r\in A^{\Gamma}$ such that $c_r=\partial b_r$, for every $r\in\mathbb R$, or, equivalently, \begin{equation}\label{b_t}\text{$c_r(g,h)=b_r(g)\sigma_g(b_r(h))b_r(gh)^{-1}$, for every $g,h\in\Gamma$.} \end{equation}
For  $r,s\in\mathbb R$, define $d_{r,s}\in A^{\Gamma}$ by letting $d_{r,s}=b_rb_sb_{r+s}^{-1}$. Since the map $\mathbb R\ni r\mapsto c_r(g,h)\in A$ is a homomorphism, for every $g,h\in\Gamma$, \eqref{b_t} implies that $d_{r,s}\in \text{Z}^1(\Gamma, A)$, for every $r,s\in\mathbb R$.  Thus, we have a Borel map  $d:\mathbb R\times\mathbb R\rightarrow\text{Z}^1(\Gamma,A)$ which is clearly a 2-cocycle. 

Next, by the hypothesis $\text{H}^1(\Gamma, A)$ is a countable discrete abelian group. 
Let $\delta:\text{Z}^1(\Gamma,A)\rightarrow\text{H}^1(\Gamma, A)$ be the quotient homomorphism. 
Since $\delta$ is continuous and $d$ is a symmetric 2-cocycle, $\delta\circ d:\mathbb R\times\mathbb R\rightarrow \text{H}^1(\Gamma,A)$ is a symmetric $2$-cocycle. 
Since $\text{H}^1(\Gamma,A)$ is countable, Lemma \ref{H^2} (1) implies that  $\delta\circ d$ is a 2-coboundary. Thus, there is a Borel map $e:\mathbb R\rightarrow \text{H}^1(\Gamma,A)$ such that \begin{equation}\label{delta}\text{$\delta(d_{r,s})=e_re_se_{r+s}^{-1}$, for $m_{\mathbb R}^2$-almost every $(r,s)\in\mathbb R^2$.}\end{equation} 

Since $\delta:\text{Z}^1(\Gamma,A)\rightarrow\text{H}^1(\Gamma,A)$ is a continuous onto homomorphism,  \cite[Theorem 2.6]{Ke94} gives a Borel map $\zeta:\text{H}^1(\Gamma,A)\rightarrow\text{Z}^1(\Gamma,A)$ such that $\delta(\zeta(e))=e$, for every $e\in \text{H}^1(\Gamma,A)$.
Thus, $f:=\zeta\circ e:\mathbb R\rightarrow\text{Z}^1(\Gamma,A)$ is a Borel map such that $e_s=\delta(f_s)$, for every $s\in\mathbb R$.
Define $k:\mathbb R\times\mathbb R\rightarrow \text{Z}^1(\Gamma,A)$ by letting $k_{r,s}=d_{r,s}(f_rf_sf_{r+s}^{-1})^{-1}$. Then $k$ is a 2-cocycle and \eqref{delta} implies that  \begin{equation}\label{k} \text{$k_{r,s}\in \text{B}^1(\Gamma,A)$, for $m_{\mathbb R}^2$-almost every $(r,s)\in\mathbb R^2$}.\end{equation}
By using  Lemma \ref{borel}, we can find a Borel map $a:\mathbb R\times\mathbb R\rightarrow A$ such that \begin{equation}\label{k}
\text{$k_{r,s}(g)=a_{r,s}^{-1}\sigma_g(a_{r,s})$, for every $g\in\Gamma$,
for $m_{\mathbb R}^2$-almost every $(r,s)\in\mathbb R^2$.} \end{equation}
Moreover, since $k_{r,s}=k_{s,r}$, we may assume that $a_{r,s}=a_{s,r}$, for every $r,s\in\mathbb R$.

View $\mathbb T<A=\text{L}^0(X,\mathbb T)$ as the closed subgroup of constant functions.
Since $k:\mathbb R\times\mathbb R\rightarrow \text{Z}^1(\Gamma,A)$ is a 2-cocycle, by combining \eqref{k} with the ergodicity of $\sigma$ it follows that \begin{equation}\text{$\mathbb T\cdot a_{r,s}a_{r+s,t}=\mathbb T\cdot a_{r,s+t}a_{s,t}$, for $m_{\mathbb R}^3$-almost every $(r,s,t)\in\mathbb R^3$.}
\end{equation}
Thus, the map $\mathbb R\times\mathbb R\ni (r,s)\mapsto \mathbb T\cdot a_{r,s}\in A/\mathbb T$ is a symmetric 2-cocycle.
By Lemma \ref{H^2} (3), we can find a Borel map $l:\mathbb R\rightarrow A$ such that \begin{equation} \text{$\mathbb T\cdot a_{r,s}=\mathbb T\cdot l_rl_sl_{r+s}^{-1}$, for almost every $(r,s)\in\mathbb R^2$.}\end{equation} Together with \eqref{k} this gives that \begin{equation}\label{l}\text{$k_{r,s}(g)=(l_r^{-1}\sigma_g(l_r))(l_s^{-1}\sigma_g(l_s))(l_{r+s}^{-1}\sigma_g(l_{r+s}))^{-1}$, for every $g\in\Gamma$ and almost every $(r,s)\in\mathbb R^2$.}
\end{equation}
Recall that $k_{r,s}=d_{r,s}(f_rf_sf_{r+s}^{-1})^{-1}=(b_rf_r^{-1})^{-1}(b_sf_s^{-1})(b_{r+s}f_{r+s}^{-1})^{-1}$, for every $r,s\in\mathbb R$. Together with \eqref{l}, we get that if $g\in\Gamma$, then the Borel map $\varphi_g:\mathbb R\rightarrow A$ given by $\varphi_g(r)=b_rf_r^{-1}l_r\sigma_g(l_r)^{-1}$ satisfies $\varphi_g(r+s)=\varphi_g(r)\varphi_g(s)$, for almost every $r,s\in\mathbb R$. Hence there is $\xi_g\in\text{L}^0(X,\mathbb R)$ such that $\varphi_g(r)=\exp(ir\xi_g)$, for almost every $r\in\mathbb R$. Thus, we derive that 
\begin{equation}\label{b_r}\text{$b_r(g)=f_r(g)l_r^{-1}\sigma_g(l_r)\exp(ir\xi_g)$, for every $g\in\Gamma$ and almost every $r\in \mathbb R$.}
\end{equation}
Since $f_r\in\text{Z}^1(\Gamma,A)$, for every $r\in\mathbb R$, from the combination  of \eqref{b_t} and \eqref{b_r} it follows that $c_r(g,h)=\exp(ir(\xi_g+\sigma_g(\xi_h)-\xi_{gh}))$, for all $g,h\in\Gamma$, and almost every $r\in\mathbb R$. This implies that $c(g,h)=\xi_g+\sigma_g(\xi_h)-\xi_{gh}$, for every $g,h\in\Gamma$, and thus $c\in \text{B}^2(\Gamma,\text{L}^0(X,\mathbb R))$.
\hfill$\blacksquare$

\begin{remark}
	Let $K$ be a locally compact abelian group. Then it is easy to see that the proof of Lemma~\ref{H^2} can in fact be applied verbatim for $K$ instead of $\R$, with the added assumption that $K$ is connected for part (1). Similarly, in Theorem~\ref{untwist}, we can consider $K$-valued cocycles $c\in$ Z$^2(\Gamma, \emph{L}^0(X, K))$, and replace $c_r(g,h)(x)=\exp(irc(g,h)(x))$ by $c_\chi(g,h)(x) = \chi(c(g,h)(x))$ for elements $\chi$ in the Pontryagin dual $\hat K$ of $K$. Assuming that $\hat K$ is connected and using the more general statement of Lemma~\ref{H^2}, the above proof of Theorem~\ref{untwist} applies verbatim to reach the same conclusion for these $K$-valued 2-cocycles.
	In particular, Theorem~\ref{untwist} holds for $K$-valued 2-cocycles whenever $K$ is of the form $\R^n\times \Lambda$ for some $n\in\N$ and some countable discrete torsion-free abelian group $\Lambda$. Indeed, these are exactly the locally compact abelian groups with connected Pontryagin dual, see for instance \cite[Theorem~26 and Corollary~4]{Mo77}.
\end{remark}

\subsection{Proof of Corollary \ref{E}}  As $\Gamma$ has property (T),   \cite[Theorem 3.4]{Sc81} implies that  $\text{B}^1(\Gamma, \text{L}^0(X,\mathbb T))$ is an open subgroup of $\text{Z}^1(\Gamma,\text{L}^0(X,\mathbb T))$. Let $c\in\text{Z}^2(\Gamma,\text{L}^0(X,\mathbb R))\setminus \text{B}^2(\Gamma,\text{L}^0(X,\mathbb R))$. For $r\in\mathbb R$, define $c_r=\exp(irc)\in\text{Z}^2(\Gamma,\text{L}^0(X,\mathbb T))$. Note that $S=\{r\in\mathbb R\mid c_r\in\text{B}^2(\Gamma,\text{L}^0(X,\mathbb T))\}$ is a subgroup of $\mathbb R$. By Theorem \ref{untwist}, $S$ is a proper subgroup of $\mathbb R$ and thus it does not contain a neighborhood of $0\in\mathbb R$. Hence, we can find a sequence $t_n\in \mathbb R\setminus S$ such that $t_n\rightarrow 0$. Since $\|c_{t_n}(g,h)-1\|_2\rightarrow 0$, for every $g,h\in\Gamma$, we can apply Theorem \ref{C} to conclude that $C^*(\Gamma)$ does not have the LP. \hfill$\blacksquare$



\subsection{Lattices with non-zero $\text{H}^2(\cdot,\mathbb R)$.} The  following lemma is needed to justify Example~\ref{example}(i).

\begin{lemma}[\!\!{\cite{BdHV08}}]\label{BdHV}
	Let $G$ be a simple Lie group with an infinite cyclic fundamental group and with property (T). If $\Gamma$ is a lattice in $G$, then $\emph{H}^2(\Gamma,\mathbb R)\not=\{0\}$.
\end{lemma}

{\it Proof.}
The proof of \cite[Corollary 3.5.6]{BdHV08} shows that there is a central extension $$0\rightarrow\mathbb Z\rightarrow\widetilde{\Gamma}\rightarrow\Gamma\rightarrow 0$$ such that $\widetilde{\Gamma}$ has property (T) and contains $\Z$ in its center. Choosing a section $q:\Gamma\rightarrow\widetilde{\Gamma}$ of the quotient homomorphism $p:\widetilde{\Gamma}\rightarrow\Gamma$, i.e. a map such that $p(q(g))=g$ for every $g\in\Gamma$, this is used to deduce that the 2-cocycle $c\in\text{Z}^2(\Gamma,\mathbb Z)$ given by $c(g,h)=q(g)q(h)q(gh)^{-1}$ does not lie in $\text{B}^2(\Gamma,\mathbb Z)$. 

To prove the conclusion, it suffices to argue that $c\notin\text{B}^2(\Gamma,\mathbb R)$. Otherwise, we can find a map $b:\Gamma\rightarrow\mathbb R$ such that $c(g,h)=b(g)b(h)b(gh)^{-1}$, for every $g,h\in\Gamma$ (here we use the multiplicative notation for addition on $\mathbb R$). 
Define $r:\widetilde{\Gamma}\rightarrow\mathbb Z$ by letting $r(g)=q(p(g))g^{-1}$, for every $g\in\widetilde{\Gamma}$.  
Since $\mathbb Z$ is central in $\widetilde{\Gamma}$, we get that $r(g)r(h)=c(p(g),p(h))r(gh)$, for every $g,h\in\widetilde{\Gamma}$. Thus, $s:\widetilde{\Gamma}\rightarrow\mathbb R$ given by $s(g)=b(p(g))^{-1}r(g)$ is a homomorphism. Since $\widetilde{\Gamma}$ has property (T), \cite[Corollary 1.3.5]{BdHV08} implies that $s(\widetilde{\Gamma})=\{0\}$. In particular, $0=s(g)=b(p(g))^{-1}q(p(g))g^{-1}=b(e)^{-1}q(e)g^{-1}$, thus $g=b(e)^{-1}q(e)$, for every $g\in\mathbb Z$. This is a contradiction and finishes the proof.
\hfill$\blacksquare$

\section{Proof of Proposition \ref{ME}} 

\subsection{An embedding result for second cohomology}
The proof of Proposition \ref{ME} relies on the following result, which appears to be of independent interest.

\begin{lemma}\label{embed}
Let $\Gamma$ be a countable group with property (T). Let $\Gamma\curvearrowright (X,\mu)$ and $\Gamma\curvearrowright (Y,\nu)$ be ergodic p.m.p. actions.  Let $c\in \emph{Z}^2(\Gamma,\emph{L}^0(X,\mathbb R))$
 and view $c\in\emph{Z}^2(\Gamma,\emph{L}^0(X\times Y,\mathbb R))$. 
 
 If $c\in\emph{B}^2(\Gamma,\emph{L}^0(X\times Y,\mathbb R))$, then $c\in \emph{B}^2(\Gamma,\emph{L}^0(X,\mathbb R))$. 
 
Thus, the natural homomorphism $\emph{H}^2(\Gamma,\emph{L}^0(X,\mathbb R))\rightarrow\emph{H}^2(\Gamma,\emph{L}^0(X\times Y,\mathbb R))$ is injective. In particular, the natural homomorphism $\emph{H}^2(\Gamma,\mathbb R)\rightarrow\emph{H}^2(\Gamma,\emph{L}^0(Y,\mathbb R))$ is injective.
\end{lemma}

To prove Lemma \ref{embed} we will need the following technical result.

\begin{lemma}\label{cohom}
Let $\Gamma$ be a countable group. Let $\Gamma\curvearrowright (X,\mu)$ and $\Gamma\curvearrowright (Y,\nu)$ be ergodic p.m.p actions. Let $\xi:X\times Y^2\rightarrow\mathbb R$ be measurable such that $\xi(x,y_1,y_2)=-\xi(x,y_2,y_1)$, for all 
$(x,y_1,y_2)\in X\times Y^2$, and suppose $\eta:X\times Y^3\rightarrow\mathbb R$ given by $\eta(x,y_1,y_2,y_3)=\xi(x,y_1,y_2)+\xi(x,y_2,y_3)+\xi(x,y_3,y_1)$ is $\Gamma$-invariant.

Then there are a $\Gamma$-invariant measurable function $\alpha:X\times Y^2\rightarrow\mathbb R$ and a measurable function $\beta:X\times Y\rightarrow\mathbb R$ with $\xi(x,y_1,y_2)=\alpha(x,y_1,y_2)+\beta(x,y_1)-\beta(x,y_2)$, for almost all $(x,y_1,y_2)\in X\times Y^2$.
\end{lemma}
Let us prove Lemma \ref{embed} assuming for the moment Lemma \ref{cohom}.

{\it Proof of Lemma \ref{embed}.} Assume that there exists a map $b:\Gamma\rightarrow\text{L}^0(X\times Y,\mathbb R)$ such that 
\begin{equation}\label{c}c(g,h)(x)=b(g)(x,y)+b(h)(g^{-1}x,g^{-1}y)-b(gh)(x,y),\end{equation} for all $g,h\in\Gamma$ and almost every $(x,y)\in X\times Y$.
Define $d:\Gamma\rightarrow\text{L}^0(X\times Y^2,\mathbb R)$ by letting $d(g)(x,y_1,y_2)=b(g)(x,y_1)-b(g)(x,y_2)$. Then \eqref{c} implies that $d\in\text{Z}^1(\Gamma,\text{L}^0(X\times Y^2,\mathbb R))$.

Since $\Gamma$ has property (T), a result of Schmidt (see \cite[Theorem 3.4]{Sc81}) and Zimmer (see \cite[Theorem 9.1.1]{Zi84}) shows that $\text{H}^1(\Gamma,\text{L}^0(Z,\mathbb R))=0$, for any ergodic p.m.p. action $\Gamma\curvearrowright (Z,\rho)$. Moreover, \cite[Proposition 2.5]{Ki17} implies that the same holds for any p.m.p. but not necessarily ergodic action $\Gamma\curvearrowright (Z,\rho)$.
Hence, there is $\xi\in\text{L}^0(X\times Y^2,\mathbb R)$ such that
\begin{equation} \label{xi}b(g)(x,y_1)-b(g)(x,y_2)=\xi(g^{-1}x,g^{-1}y_1,g^{-1}y_2)-\xi(x,y_1,y_2), 
\end{equation}
for all $g\in\Gamma$ and almost every $(x,y_1,y_2)\in X\times Y^2$. Since \eqref{xi} still holds if we replace $\xi$ with the function $X\times Y^2\ni (x,y_1,y_2)\mapsto \frac{\xi(x,y_1,y_2)-\xi(x,y_2,y_1)}{2}\in\mathbb R$, we may additionally assume that 
\begin{equation}\label{anti}\text{$\xi(x,y_1,y_2)=-\xi(x,y_2,y_1)$, for all $(x,y_1,y_2)\in X\times Y^2$.}
\end{equation}
Next, by \eqref{xi}, $\eta:X\times Y^3\rightarrow \mathbb R$ given by $\eta(x,y_1,y_2,y_3)=\xi(x,y_1,y_2)+\xi(x,y_2,y_3)+\xi(x,y_3,y_1)$ is $\Gamma$-invariant. By applying Lemma \ref{cohom}, we can find a $\Gamma$-invariant measurable function $\alpha:X\times Y^2\rightarrow\mathbb R$ and a measurable function $\beta:X\times Y\rightarrow\mathbb R$ such that $\xi(x,y_1,y_2)=\alpha(x,y_1,y_2)+\beta(x,y_1)-\beta(x,y_2)$, for almost all $(x,y_1,y_2)\in X\times Y^2$.
 In combination with \eqref{xi}, it follows that for every $g\in\Gamma$ the map $X\times Y\ni (x,y)\mapsto b(g)(x,y)-(\beta(g^{-1}x,g^{-1}y)-\beta(x,y))\in\mathbb R$ 
is independent of the $y$-variable. Hence, there is $a(g)\in \text{L}^0(X,\mathbb R)$ such that $b(g)(x,y)-(\beta(g^{-1}x,g^{-1}y)-\beta(x,y))=a(g)(x)$, for almost every $(x,y)\in X\times Y$. Together with \eqref{c}, we get that $c(g,h)(x)=a(g)(x)+a(h)(g^{-1}x)-a(gh)(x)$, for all $g,h\in\Gamma$ and almost every $x\in  X$. This proves that $c\in\text{B}^2(\Gamma,\text{L}^0(X,\mathbb R))$, as desired.
\hfill$\blacksquare$

\begin{remark}
Lemma \ref{cohom} has a short proof if $\Gamma\curvearrowright (Y,\nu)$ is weakly mixing. In this case, the product action $\Gamma\curvearrowright (X\times Y^3, \mu\times\nu^3)$ is ergodic and so $\eta$ is constant.  Since $\eta(x,y_1,y_3,y_2)=-\eta(x,y_1,y_2,y_3)$, we get that $\eta(x,y_1,y_2,y_3)=0$, for almost every $(x,y_1,y_2,y_3)\in X\times Y^3$. By Fubini's theorem, there is $z\in Y$ such that $\eta(x,y_1,y_2,z)=0$ and thus $\xi(x,y_1,y_2)=\xi(x,y_1,z)-\xi(x,y_2,z)$, for almost every $(x,y_1,y_2)\in X\times Y^2$. To prove Lemma \ref{cohom} in general, we will need a much more involved argument, working with the maximal compact factors of the actions considered.

\end{remark}

{\it Proof of Lemma \ref{cohom}.} First, we argue that the actions $\Gamma\curvearrowright (X,\mu)$ and $\Gamma\curvearrowright (Y,\nu)$ may be assumed compact. Let $\Gamma\curvearrowright (X_0,\mu_0)$ and $\Gamma\curvearrowright (Y_0,\nu_0)$ be their maximal compact factors. Let $p:X\rightarrow X_0$ and $q:Y\rightarrow Y_0$ be the  $\Gamma$-equivariant factor maps. Since the Koopman representations of $\Gamma$ on $\text{L}^2(X)\ominus \text{L}^2(X_0)$ and $\text{L}^2(Y)\ominus \text{L}^2(Y_0)$ are weakly mixing, any $\Gamma$-invariant map $\eta:X\times Y^3\rightarrow\mathbb R$ factors through $p\times q^3$. Let $\rho:X_0\times Y_0^3\rightarrow\mathbb R$ be measurable such that $\eta=\rho\circ(p\times q^3)$. By Fubini's theorem, we can find $z\in Y$ such that $\rho(p(x),q(y_1),q(y_2),q(z))=\xi(x,y_1,y_2)+\xi(x,y_2,z)-\xi(x,y_1,z)$, for almost every $(x,y_1,y_2)\in X\times Y^2$. Then $\xi_0:X\times Y^2\rightarrow \mathbb R$ given by $\xi_0(x,y_1,y_2)=\rho(p(x),q(y_1),q(y_2),q(z))$ factors through $p\times q^2$ and satisfies the same properties as $\xi$. Since it suffices to prove the conclusion for $\xi_0$ instead of $\xi$, we may indeed assume that $X=X_0$ and $Y=Y_0$.

Since the ergodic actions $\Gamma\curvearrowright (X,\mu)$ and $\Gamma\curvearrowright (Y,\nu)$ are  compact, they are isomorphic to left translation actions $\Gamma\curvearrowright G/K$ and $\Gamma\curvearrowright G/L$, where  $G$ is a compact group containing $\Gamma$ densely and $K,L<G$ are closed subgroups. As $\eta$ is $\Gamma$-invariant and $\Gamma<G$ is dense,  $\eta$ is $G$-invariant. Then  Fubini's theorem gives $a\in G$ such that $\eta(gaK,gy_1,gy_2,gy_3)=\eta(aK,y_1,y_2,y_3)$, for almost every $(g,y_1,y_2,y_3)\in G\times Y^3$.
By applying Fubini's theorem again we find $b\in G$ such that $$\eta(bK,gy_1,gy_2,gy_3)=\eta(gaK,gy_1,gy_2,gy_3)=\eta(aK,y_1,y_2,y_3),$$ for almost every $(g,y_1,y_2,y_3)\in bKa^{-1}\times Y^3$. Thus, $\eta(bK,gy_1,gy_2,gy_3)=\eta(bK,hy_1,hy_2,hy_3)$, for almost every $(g,h,y_1,y_2,y_3)\in (bKa^{-1})^2\times Y^3$. 
Hence, if we let $K_0=bKb^{-1}$, then  \begin{equation}\label{eta}\text{$\eta(bK,ky_1,ky_2,ky_3)=\eta(bK,y_1,y_2,y_3)$, for almost every $(k,y_1,y_2,y_3)\in K_0\times Y^3$.}\end{equation} 
Define $\rho\in\text{Z}^1(K_0,\text{L}^0(Y^2,\mathbb R))$ by $\rho(k)(y_1,y_2)=\xi(bK,k^{-1}y_1,k^{-1}y_2)-\xi(bK,y_1,y_2)$. Then \eqref{eta} implies that $\rho(k)(y_1,y_2)+\rho(k)(y_2,y_3)+\rho(k)(y_3,y_1)=0$, for almost every $(k,y_1,y_2,y_3)\in K_0\times Y^3$. Moreover, we have $\rho(k)(y_1,y_2)=-\rho(k)(y_2,y_1)$, for every $(k,y_1,y_2)\in K_0\times Y^2$.
By Fubini's theorem we can thus find $y_3\in Y$ such that the map $\sigma:K_0\rightarrow\text{L}^0(Y,\mathbb R)$ given by $\sigma(k)(y)=\rho(k)(y,y_3)$ satisfies
\begin{equation}\label{sigma}
\text{$\sigma(k)(y_1)-\sigma(k)(y_2)=\rho(k)(y_1,y_2)$, for almost every $(k,y_1,y_2)\in K_0\times Y^2$.}
\end{equation}
Since  $\rho\in\text{Z}^1(K_0,\text{L}^0(Y^2,\mathbb R))$,  \eqref{sigma} implies that for almost every $(k,k',y_1,y_2)\in K_0^2\times Y^2$ we have $$\sigma(k)(y_1)+\sigma(k')(k^{-1}y_1)-\sigma(kk')(y_1)=\sigma(k)(y_2)+\sigma(k')(k^{-1}y_2)-\sigma(kk')(y_2).$$ Hence, for almost every $(k,k')\in K_0^2$, there is $c(k,k')\in\mathbb R$ such that  \begin{equation}\label{sigma2}\text{$\sigma(k)(y)+\sigma(k')(k^{-1}y)-\sigma(kk')(y)=c(k,k')$, for almost every $y\in Y$.}\end{equation}
Then one checks that $c\in\text{Z}^2(K_0,\mathbb R)$. Since $K_0$ is a compact group, \cite[Theorem A]{AM13} implies that $\text{H}^2(K_0,\mathbb R)=0$.  Thus, we can find  $b:K_0\rightarrow\mathbb R$ measurable such that $c(k,k')=b(k)+b(k')-b(kk')$, for almost every $(k,k')\in K_0^2$. In combination with \eqref{sigma2} it follows that the map $\tau:K_0\rightarrow\text{L}^0(Y,\mathbb R)$ given by $\tau(k)(y)=\sigma(k)(y)-b(k)$ belongs to $\text{Z}^1(K_0,\text{L}^0(Y,\mathbb R))$. Using again that $K_0$ is  a compact group and \cite[Theorem A]{AM13} we get that $\text{H}^1(K_0,\text{L}^0(Y,\mathbb R))=0$. Thus, there is $f\in\text{L}^0(Y,\mathbb R)$ such that $\tau(k)(y)=f(k^{-1}y)-f(y)$, for almost every $(k,y)\in K_0\times Y$. By combining this with \eqref{sigma} and the definition of $\rho$, we get that for almost every $(k,y_1,y_2)\in K_0\times Y^2$ we have$$\xi(bK,k^{-1}y_1,k^{-1}y_2)-\xi(bK,y_1,y_2)=(f(k^{-1}y_1)-f(y_1))-(f(k^{-1}y_2)-f(y_2)).$$ Thus, the map $\zeta:Y^2\rightarrow \mathbb R$ given by $\zeta(y_1,y_2)=\xi(bK,y_1,y_2)-f(y_1)+f(y_2)$ is $K_0$-invariant. Using the fact that $\eta$ is $G$-invariant, for almost every $(x,y_1,y_2,y_3)\in G\times Y^3$, we get that 
\begin{align*}&\xi(xK,y_1,y_2)+\xi(xK,y_2,y_3)+\xi(xK,y_3,y_1)\\&=\xi(bK,bx^{-1}y_1,bx^{-1}y_2)+\xi(bK,bx^{-1}y_2,bx^{-1}y_3)+\xi(bK,bx^{-1}y_3,bx^{-1}y_1)\\&=\zeta(bx^{-1}y_1,bx^{-1}y_2)+\zeta(bx^{-1}y_2,bx^{-1}y_3)+\zeta(bx^{-1}y_3,bx^{-1}y_1).
\end{align*}
Note that the map $G\times Y\ni (x,z)\mapsto (x,xz)\in G\times Y$ is measure preserving. 
By applying Fubini's theorem, we can find $z\in Y$ such that for almost every $(x,y_1,y_2)\in G\times Y^2$, the last displayed identity holds for $y_3=xz$. Equivalently, for almost every $(x,y_1,y_2)\in G\times Y^2$ we have that $$\xi(xK,y_1,y_2)+\xi(xK,y_2,xz)-\xi(xK,y_1,xz)=\zeta(bx^{-1}y_1,bx^{-1}y_2)+\zeta(bx^{-1}y_2,bz)-\zeta(bx^{-1}y_1,bz).$$
Thus, if we define $F:G\times Y\rightarrow\mathbb R$ by letting $F(x,y)=\xi(xK,y,xz)-\zeta(bx^{-1}y,bz)$ then we have \begin{equation}\label{zeta}\xi(xK,y_1,y_2)-\zeta(bx^{-1}y_1,bx^{-1}y_2)=F(x,y_1)-F(x,y_2),\end{equation} for almost every $(x,y_1,y_2)\in G\times Y^2$. 

Since $\zeta$ is $K_0$-invariant, the map $G\times Y^2\ni (x,y_1,y_2)\mapsto \zeta(bx^{-1}y_1,bx^{-1}y_2)\in\mathbb R$ factors through $G\times Y^2\rightarrow G/K\times Y^2$. Define $\alpha:G/K\times Y^2\rightarrow\mathbb R$ by letting $\alpha(xK,y_1,y_2)=\zeta(bx^{-1}y_1,bx^{-1}y_2)$, and note that $\alpha$ is $\Gamma$-invariant.

 By combining this fact with \eqref{zeta} we get that $F(xk,y_1)-F(xk,y_2)=F(x,y_1)-F(x,y_2)$, for almost every $(x,y_1,y_2,k)\in G\times Y^2\times K$. Hence, there is $w:K\times G\rightarrow\mathbb R$ measurable such that $F(xk,y)-F(x,y)=w(k,x)$, for almost every $(x,y,k)\in G\times Y\times K$. Then $w$ is a 1-cocycle for the right translation action of $K$ on $G$. Since $K$ is compact, \cite[Theorem A]{AM13} gives a measurable map $v:G\rightarrow\mathbb R$ such that $w(k,x)=v(xk)-v(x)$, for almost every $(x,k)\in G\times K$. Thus, the map $G\times Y\ni (x,y)\mapsto F(x,y)-v(x)\in\mathbb R$ factors through $G\times Y\rightarrow G/K\times Y$. Let $\beta:G/K\times Y\rightarrow\mathbb R$ be given by $\beta(xK,y)=F(x,y)-v(x)$. Finally, \eqref{zeta} rewrites as 
$\xi(xK,y_1,y_2)-\alpha(xK,y_1,y_2)=\beta(xK,y_1)-\beta(xK,y_2)$, for almost every $(x,y_1,y_2)\in G\times Y^2$. This finishes the proof. 
\hfill$\blacksquare$

\subsection{Proof of Proposition \ref{ME}} We begin by introducing some notation and recalling several facts that we will use in the proof.
For a countable p.m.p. equivalence relation $\mathcal R$ on a probability space $(T,\lambda)$ we denote by $\text{H}^2(\mathcal R,\mathbb R)$ the second cohomology  group of $\mathcal R$ with values in $\mathbb R$ \cite{FM77} and by $\mathcal R_{|T_0}=\mathcal R\cap (T_0\times T_0)$ its restriction to a measurable set $T_0\subset T$.
If $T_0\subset T$ is a measurable set and $T_1$ is the $\mathcal R$-saturation of $T_0$ (i.e., the set of $x\in T$ such that there is $y\in T_0$ with $(x,y)\in\mathcal R$), then \cite[Lemma 2.2]{Ki17} shows that $\text{H}^2(\mathcal R_{|T_0},\mathbb R)\cong \text{H}^2(\mathcal R_{|T_1},\mathbb R)$.
Finally, if $\Sigma\curvearrowright (T,\lambda)$ is a p.m.p. action of a countable group, we denote by $\mathcal R(\Sigma\curvearrowright T)=\{(x,y)\in T\times T\mid \Sigma\cdot x=\Sigma\cdot y\}$ its orbit equivalence relation. By \cite[Theorem 5]{FM77}, if $\Sigma$ acts freely, then $\text{H}^2(\mathcal R(\Sigma\curvearrowright T),\mathbb R)\cong \text{H}^2(\Sigma,\text{L}^0(T,\mathbb R))$.

Let $\Gamma$ be a countable group with property (T) admitting an ergodic p.m.p. action $\Gamma\curvearrowright (X,\mu)$ with $\text{H}^2(\Gamma,\text{L}^0(X,\mathbb R))\not =0$. Let $\Lambda$ be a countable group that is measure equivalent to $\Gamma$.  
By \cite[Corollary 1.4]{Fu99a}, $\Lambda$ has property (T). Moreover, by \cite[Lemma 3.2 and Theorem 3.3]{Fu99b}, there exist free ergodic p.m.p. actions $\Gamma\curvearrowright (Y,\nu)$, $\Lambda\curvearrowright (Z,\rho)$ that are stably orbit equivalent: there are non-negligible measurable sets $Y_0\subset Y$, $Z_0\subset Z$ and a measure preserving isomorphism $\theta:Y_0\rightarrow Z_0$ 
such that  $(\theta\times\theta)(\mathcal R(\Gamma\curvearrowright Y)_{|Y_0})=\mathcal R(\Lambda\curvearrowright Z)_{|Z_0}$. 

Consider the product action $\Gamma\curvearrowright (X\times Y,\mu\times\nu)$. By using the previous paragraph, one can construct a free p.m.p. action $\Lambda\curvearrowright (W,\eta)$ and a non-negligible measurable set $W_0\subset W$ such that $\mathcal R(\Gamma\curvearrowright X\times Y)_{|(X\times Y_0)}$ is isomorphic to $\mathcal R(\Lambda\curvearrowright W)_{|W_0}$ (see \cite[Claim 5, Section 3]{Io11}). Since the action $\Gamma\curvearrowright (Y,\nu)$ is ergodic, we have that $\Gamma\cdot (X\times Y_0)=X\times Y$. Denote $W_1=\Lambda\cdot W_0$. Then $W_1$ is $\Lambda$-invariant and
by combining the facts from the first paragraph of the proof we get that 
\begin{align*}\text{H}^2(\mathcal R(\Lambda\curvearrowright W_1),\mathbb R)&\cong \text{H}^2(\mathcal R(\Lambda\curvearrowright W)_{|W_0},\mathbb R)\\&\cong \text{H}^2(\mathcal R(\Gamma\curvearrowright X\times Y)_{|X\times Y_0},\mathbb R)\\&\cong \text{H}^2(\mathcal R(\Gamma\curvearrowright X\times Y),\mathbb R)\\&\cong\text{H}^2(\Gamma,\text{L}^0(X\times Y,\mathbb R)).\end{align*}
Since $\text{H}^2(\Gamma,\text{L}^0(X,\mathbb R))\not =0$, using Lemma \ref{embed} we conclude that $\text{H}^2(\mathcal R(\Lambda\curvearrowright W_1),\mathbb R)\not=0$. By using the ergodic decomposition  of the free p.m.p. action $\Lambda\curvearrowright (W_1,\eta(W_1)^{-1}\eta_{|W_1})$ and \cite[Proposition 2.5]{Ki17} we can find a $\Lambda$-invariant probability measure $\eta_1$ on $W_1$ such the p.m.p. action $\Lambda\curvearrowright (W_1,\eta_1)$ is free and ergodic, and its orbit equivalence relation $\mathcal S$ satisfies $\text{H}^2(\mathcal S,\mathbb R)\not=0$. Finally, we have that $\text{H}^2(\Lambda,\text{L}^0((W_1,\eta_1),\mathbb R))\cong\text{H}^2(\mathcal S,\mathbb R)\not=0$, which proves the conclusion.
\hfill$\blacksquare$

\subsection{An embedding result for first cohomology}

\begin{theorem}\label{H1injective}
Let $\Gamma\curvearrowright^{\sigma} (X,\mu)$ be a p.m.p. action of a countable group $\Gamma$ with spectral gap. 

Then the natural homomorphism  $\emph{H}^1(\Gamma,\emph{L}^2(X,\mathbb R))\rightarrow \emph{H}^1(\Gamma,\emph{L}^0(X,\mathbb R))$ is injective.

\end{theorem}

{\it Proof.}  Let $c:\Gamma\rightarrow\text{L}^2(X,\mathbb R)$ be a 1-cocycle for which there is $b\in\text{L}^0(X,\mathbb R)$ such that $c(g)=\sigma_g(b)-b$, for all $g\in\Gamma$.  Let $M>0$ such that $Y=\{x\in X\mid |b(x)|<M\}$ satisfies $\mu(Y)>\frac{1}{2}$. 
Since $\sigma$ has spectral gap, we can find $F\subset\Gamma$ finite and $C_1>0$ such that \begin{equation}\label{spectral}
\text{$\|\xi-\int_X\xi\;\text{d}\mu\cdot 1\|_2\leq C_1\max_{g\in F}\|\sigma_g(\xi)-\xi\|_2$, for every $\xi\in\text{L}^2(X).$}\end{equation}

For every $t\in\mathbb R$, let $u_t=\exp(itb)\in\text{L}^0(X,\mathbb T)$ and $\alpha_t=\int_Xu_t\;\text{d}\mu$. Then $\sigma_g(u_t)u_t^*=\exp(itc(g))$.
Using that $|\exp(ix)-1|\leq |x|$, for all $x\in\mathbb R$, $\|\sigma_g(u_t)-u_t\|_2=\|\sigma_g(u_t)u_t^*-1\|_2$ and \eqref{spectral}  we get that
$$
\|u_t-\alpha_t\cdot 1\|_2\leq C_1\max_{g\in F}\|\sigma_g(u_t)-u_t\|_2=C_1 \|\exp(itc(g))-1\|_2\leq C_1t\max_{g\in F}\|c(g)\|_2.
$$
Thus, if we let $C_2=C_1\max_{g\in F}\|c(g)\|_2<\infty$, then \begin{equation}\label{u_t}\text{$\|u_t-\alpha_t\cdot 1\|_2\leq C_2t$, for every $t\in\mathbb R$.}\end{equation}If $t\in \mathbb R$, then since $\int_Y|\exp(itb(x))-\alpha_t|^2\;\text{d}\mu(x)\leq \|u_t-\alpha_t\cdot 1\|_2^2\leq C_2^2t^2,$ there is $x\in Y$ such that $|\exp(itb(x)-\alpha_t|\leq \sqrt{2}C_2|t|$. Since $x\in Y$, we have that $|\exp(itb(x))-1|\leq |t b(x)|\leq M|t|$. Thus, we get  $|\alpha_t-1|\leq C_3t$, where $C_3=\sqrt{2}C_2+M$. Letting $C_4=C_2+C_3$ and using \eqref{u_t} this gives that $\|u_t-1\|_2\leq C_4t$, for every $t\in\mathbb R$.

Finally, for $\delta>0$, let $X_{\delta}=\{x\in X\mid |b(x)|\leq \delta\}$. Since $\int_{X_\delta}\frac{|\exp(itb(x))-1|^2}{t^2}\;\text{d}\mu(x)\leq C_4^2$, $\frac{|\exp(itb(x))-1|^2}{t^2}\leq \delta^2$ and $\lim\limits_{t\rightarrow 0}\frac{|\exp(itb(x))-1|^2}{t^2}=|b(x)|^2$, for every $t\in\mathbb R$ and $x\in X_{\delta}$, the dominated convergence theorem implies that $\int_{X_{\delta}}|b(x)|^2\;\text{d}\mu(x)\leq C_4^2$. Since this holds for every $\delta>0$, the monotone convergence theorem implies that $b\in \text{L}^2(X,\mathbb R)$, which finishes the proof.
\hfill$\blacksquare$

\begin{corollary} \emph{(cf. \cite[Corollary 4.2]{PS12})}
Let $\Gamma\curvearrowright (X,\mu)$ be an ergodic p.m.p. action that is $\{\mathbb T\}$-cocycle superrigid. 

Then $\emph{H}^1(\Gamma,\emph{L}^0(X,\mathbb R))=0$. Moreover, if $\Gamma\curvearrowright (X,\mu)$ has spectral gap, then $\emph{H}^1(\Gamma,\emph{L}^2(X))=0$.

\end{corollary}

{\it Proof.}  The main assertion is a direct consequence of a result from \cite{MS80}. Let $c\in\text{Z}^1(\Gamma,\text{L}^0(X,\mathbb R))$. For $r\in\mathbb R$, define $c_r\in\text{Z}^1(\Gamma,\text{L}^0(X,\mathbb T))$ by $c_r(g)=\exp(irc(g))$. Since  $\Gamma\curvearrowright (X,\mu)$ is $\{\mathbb T\}$-cocycle superrigid, $c_r$ is cohomologous to a character $\eta_r:\Gamma\rightarrow\mathbb T$, for every $r\in\mathbb R$. Then applying \cite[Theorem 6.2, equivalence of conditions (6) and (7)]{MS80} implies that $c\in \text{B}^1(\Gamma,\text{L}^0(X,\mathbb R))$, as desired. The moreover assertion now follows by combining the main assertion and Theorem \ref{H1injective}.
\hfill$\blacksquare$

\section{Proof of Theorem \ref{infpres}}
Let $\Gamma$ be a non-finitely presented group with property (T). Assume by contradiction that $C^*(\Gamma)$ has the \text{LP}.
As $\Gamma$ is finitely generated \cite{Ka67}, 
it has a presentation $\Gamma=\langle s_1,...,s_k\mid r_l, l\in\mathbb N\rangle$.
By a theorem of Shalom (see \cite[Theorem 6.7]{Sh00}), 
there exists $m\in\mathbb N$ such that the finitely presented group $\Gamma_0=\langle s_1,...,s_k\mid r_1,...,r_m\rangle$ has property (T).  Denote by $\pi:\Gamma_0\rightarrow\Gamma$ the canonical onto homomorphism. Let $\delta:\Gamma\rightarrow\Gamma_0$ be a map such that $\pi(\delta(g))=g$, for every $g\in\Gamma$.

For $n\in\mathbb N$, let $\Gamma_n=\langle s_1,...,s_k\mid r_1,...,r_{m+n}\rangle$ and $\pi_n:\Gamma_0\rightarrow\Gamma_n$ be  the canonical onto homomorphism. Let $M_n=\text{L}(\Gamma_n)$ and $\tau_n:M_n\rightarrow\mathbb C$ be its usual trace.
Define $\varphi_n:\Gamma\rightarrow\mathcal U(M_n)$ by letting $\varphi_n(g)=u_{\pi_n(\delta(g))}$.
If $g,h\in\Gamma$, then $\delta(gh)^{-1}\delta(g)\delta(h)\in\ker(\pi)$. Thus, $\delta(gh)^{-1}\delta(g)\delta(h)\in\ker(\pi_n)$ so $\varphi_n(gh)=\varphi_n(g)\varphi_n(h)$, for large $n$. 
Thus, $\varphi_n:\Gamma\rightarrow\mathcal U(M_n)$, $n\in\mathbb N$, is an asymptotic homomorphism.

Since $C^*(\Gamma)$ has the {\text LP}, Corollary \ref{ucplift} implies the existence of u.c.p. maps $\Phi_n:C^*(\Gamma)\rightarrow \text{L}(\Gamma_n)$, $n\in\mathbb N$, such that $\|\Phi_n(u_g)-u_{\pi_n(\delta(g))}\|_{2,\tau_n}=\|\Phi_n(u_g)-\varphi_n(g)\|_{2,\tau_n}\rightarrow 0$, for every $g\in\Gamma$. Thus, \begin{equation}\label{descend}\text{$\|\Phi_n(u_{\pi(g)})-u_{\pi_n(\delta(\pi(g)))}\|_{2,\tau_n}\rightarrow 0$, for every $g\in\Gamma_0$.}\end{equation} If $g\in\Gamma_0$, then $g^{-1}\delta(\pi(g))\in\ker(\pi)$. Thus, $g^{-1}\delta(\pi(g))\in\ker(\pi_n)$ and hence $\pi_n(\delta(\pi(g)))=\pi_n(g)$, for $n$ large enough. In combination with \eqref{descend}, we deduce that \begin{equation}\label{descend2}\text{$\|\Phi_n(u_{\pi(g)})-u_{\pi_n(g)}\|_{2,\tau_n}\rightarrow 0$, for every $g\in\Gamma_0$.}\end{equation}
For $n\in\mathbb N$, let $\psi_n:\Gamma_0\rightarrow\mathbb C$ be the positive definite map given by $\psi_n(g)=\tau_n(\Phi_n(u_{\pi(g)})u_{\pi_n(g)}^*)$. By \eqref{descend2}, $\psi_n(g)\rightarrow 1$, for any $g\in\Gamma_0$. Since $\Gamma_0$ has property (T), we get $\sup_{g\in\Gamma_0}|\psi_n(g)-1|\rightarrow 0$.
Since $\|\Phi_n(u_{\pi(g)})\|\leq 1$, we get that $\|\Phi_n(u_{\pi(g)})-u_{\pi_n(g)}\|_{2,\tau_n}^2\leq 2 |\psi_n(g)-1|$, for every $g\in\Gamma_0$. Hence, we derive that $\sup_{g\in\Gamma_0}\|\Phi_n(u_{\pi(g)})-u_{\pi_n(g)}\|_{2,\tau_n}\rightarrow 0$. As  a consequence, there is $n\in\mathbb N$ such that $\|\Phi_n(u_{\pi(g)})-u_{\pi_n(g)}\|_{2,\tau_n}<1$, for every $g\in\Gamma_0$. Thus, if $g\in\ker(\pi)$, then since $\Phi_n$ is unital we get that $\|1-u_{\pi_n(g)}\|_{2,\tau_n}<1$ and therefore $g\in\ker(\pi_n)$. This implies that $\ker(\pi)=\ker(\pi_n)$, hence $\Gamma=\Gamma_n$ would be finitely presented, which is a contradiction.
\hfill$\blacksquare$

\end{document}